\documentclass[11pt]{amsart}
\usepackage[utf8]{inputenc}
\usepackage[english]{babel}
\usepackage[active]{srcltx}
\usepackage{verbatim}
\usepackage{graphicx}
\usepackage{psfrag}
\usepackage{subfig}

\usepackage{amsmath,latexsym,amssymb,verbatim}

\newtheorem{lema}{Lemma}[section]
\newtheorem{theo}[lema]{Theorem}

\newtheorem{prop}[lema]{Proposition}
\newtheorem{coro}[lema]{Corollary}
\theoremstyle{definition}
\newtheorem*{defi}{Definition}
\theoremstyle{remark}
\newtheorem{rema}[lema]{Remark}

\newtheorem{exam}[lema]{Example}
\newtheorem{problem}[lema]{Problem}
\newtheorem{conj}[lema]{Conjecture}

\def\sideremark#1{\ifvmode\leavevmode\fi\vadjust{\vbox to0pt{\vss % the remark
      \hbox to 0pt{\hskip\hsize\hskip1em           %                will appear only
 \vbox{\hsize2cm\tiny\raggedright\pretolerance10000%                on the side
 \noindent #1\hfill}\hss}\vbox to8pt{\vfil}\vss}}}%
\newcommand{\R}{\mathbb{R}}
\newcommand{\C}{\mathbb{C}}
\newcommand{\Q}{\mathbb{Q}}

\newcommand{\id}{\text{id}}

\address{
$^1 $ 
Departamento de Matem\'aticas
Universidad de Extremadura
06006 BADAJOZ
Spain,
$^2$ Universit\'e de Bourgogne, Institut de Math\'ematiques de Bourgogne, UMR 5584 du CNRS,
UFR Sciences et Techniques, 9, av. A. Savary, BP 47870, 21078 DIJON CEDEX, France
}

\begin{document}
\title[Vanishing Abelian Integrals on Zero-Dimensional Cycles]{Vanishing Abelian Integrals on Zero-Dimensional Cycles }
\author{A. \'Alvarez$^1$, J.L. Bravo$^1$, P. Marde\v si\'c$^2$}

\begin{abstract}

In this paper we study conditions for the vanishing of Abel\-ian integrals on families of zero-dimensional cycles. That is, for any rational function $f(z)$, characterize all rational functions $g(z)$ and zero-sum integers $\{n_i\}$ such that the function $t\mapsto\sum n_ig(z_i(t))$ vanishes identically. Here $z_i(t)$ are continuously depending roots of $f(z)-t$. We introduce a notion of (un)balanced cycles. Our main result is an inductive solution of the problem of vanishing of Abelian integrals when $f,g$ are polynomials on a family of zero-dimensional cycles under the assumption that the family of cycles we consider is unbalanced as well as all the cycles encountered in the inductive process. We also solve the problem on some balanced cycles. 

The main motivation for our study  is the problem of vanishing of Abelian integrals on single families of  one-dimensional cycles. We show that our problem and our main result are sufficiently rich to include some related problems, as hyper-elliptic integrals on one-cycles, some applications to slow-fast planar systems,
and the polynomial (and trigonometric) moment problem for Abel equation. This last problem was recently solved by Pakovich and Muzychuk (\cite{PM} and \cite{P}). Our approach is largely inspired by their work, thought we provide examples of vanishing Abelian integrals on zero-cycles which are not given as a sum of composition terms contrary to the situation in the solution of the polynomial moment problem.

\end{abstract}

\subjclass[2010]{34C07, 34C08, 34D15, 34M35}

\keywords{Abelian integral, persistent center, tangential center, moment problem, center-focus problem, canard center, slow-fast system}

\thanks{
The first author was partially supported by Junta de Extremadura
and FEDER funds. The second author was partially supported by
Junta de Extremadura and a MCYT/FEDER grant number MTM2008-05460.
The first two authors are grateful to the Universit\'{e} de
Bourgogne, for the hospitality and support during the visit when
this work was started. The third author thanks the Universidad de
Extremadura for the hospitality and support during the visit when this work was finished. }

\maketitle

\section{Introduction}

Zero-dimensional Abelian integrals were introduced by Gavrilov and Mova\-sati in \cite{GM}  and  the problem of the bound for the number of their zeros was studied. In this paper we study conditions for their  identical vanishing on families of zero-dimensional cycles.

\medskip

\begin{defi}\label{abelian}
Let $f\in\mathbb{C}(z)$ be a rational function of degree $m>1$. Let $(z_i(t))_{1\leq i\leq m}$ denote an $m$-tuple of analytic preimages  $z_i(t)\in f^{-1}(t)$, where $t\in\mathbb{C}\backslash
\Sigma$ and $\Sigma$ is the set of critical values of $f$. Define a \emph{zero-dimensional cycle} (shorter \emph{cycle}) $C(t)$ of $f$ as the sum 
\begin{equation}\label{cycle}
C(t)=\sum_{i=1}^m n_i z_i(t),\quad \text{with} \quad \sum_{i=1}^m n_i=0, \quad n_i\in\mathbb{Z}.
\end{equation}
Zero-dimensional cycles form the reduced $0$-th homology group denoted $\tilde{H}_0(f^{-1}(t))$.
We say that a cycle $C(t)$ is \emph{simple} if it is of the form $C(t)=z_j(t)-z_i(t)$,
and \emph{trivial} if $n_i=0$ for every $i$.

Let $g \in \mathbb{C}(z)$ be a rational function. As in \cite{GM}, define \emph{zero-dimensional Abelian integrals of $g$ along the cycle $C(t)$} by
\begin{equation}\label{0int}
\int_{C(t)}g:=\sum_{i=1}^m n_i g(z_i(t)).
\end{equation}
\end{defi}
Note that Abelian integrals on zero-cycles are simply algebraic functions, contrary to the one-dimensional case.

\medskip

We study two problems:
\begin{problem}\label{problem1}
Characterize rational functions $f$, $g$  such that for any family of cycles $C(t)$ of $f$ the Abelian integral $\int_{C(t)}g$ vanishes identically.
\end{problem}

\begin{problem}\label{problem2}(Tangential center problem)
Characterize rational functions $f$, $g$ and cycles $C(t)$ of $f$ such that the Abelian integral $\int_{C(t)}g$ vanishes identically.
\end{problem}

In Section~\ref{sec:main} we state our results on those two problems. Concretely, we solve completely Problem~\ref{problem1} and give an inductive solution of Problem~\ref{problem2} under the assumptions that $f,g$ are polynomials, and that the family of cycles we consider is unbalanced as well as all the cycles encountered in the inductive process. The definition of balanced and unbalanced cycles is related to the behavior
of the branches $z_i(t)$ as $t$ tends to infinity, and it is made precise in Section~\ref{sec:main}. We also solve the problem on some balanced cycles. In particular we provide examples  of vanishing of Abelian integrals on zero cycles which are not given as a sum of composition terms contrary to the situation in the polynomial moment problem as solved by Pakovich and Muzychuk. The proofs of these results are included in Sections~\ref{sec:all}, \ref{sec:t2} and \ref{balanced}.

In Section~\ref{sec:3} we show that problems above are the zero-dimensional version of the analogous problems on the vanishing of Abelian integrals on one-cycles which appear in the study of periodic solutions of deformations of Hamiltonian vector fields in the plane. 

We also push further the analogy between deformations of planar integrable systems and deformations of polynomials in one variable. In particular we define the displacement function along a zero-cycle,
zero-dimensional Abelian integrals being its principal part. We study identical vanishing of the displacement function.

Problem~\ref{problem1} for one-dimensional Abelian integrals was studied by Gavrilov \cite{G} and completely solved by Bonnet-Dimca \cite{BD} in the polynomial case and by Muci\~no in \cite{Muc} in the generic rational case. Problem~\ref{problem2} for one-dimensional \emph{generic} Abelian integrals was solved by Ilyashenko \cite{I}. It was solved in the hyper-elliptic case (which is non-generic) for vanishing cycles by Christopher and the third author \cite{CM}. The general one-dimensional case remains open for both problems.

In Sections~\ref{hyper}, \ref{sec:sf} we give applications of our results to study the one-dimensional tangential center problem and related problems. Concretely, we study the vanishing of hyper-elliptic integrals on one-cycles and some related problems for slow-fast planar systems.

In Section~\ref{sec:moment} we prove that the polynomial and the trigonometric moment problem can be stated as the vanishing of zero-dimensional Abelian integrals on special cycles. In particular, in the polynomial case the cycle is unbalanced and so are all the cycles in the inductive process so our results can be applied. The moment problem was started by Briskin, Françoise and Yomdim (\cite{BFY1,BFY}) and it  has received a lot of attention (see e.g. \cite{BRY,C,PCounter,P2}), but was only recently solved by Pakovich and Muzychuk (\cite{PM} and \cite{P}). Our approach is largely inspired by their work.

In \cite{BFY1} Briskin, Fran\c coise, Yomdin formulated the composition conjecture for the moment problem (see Section \ref{sec:moment}). Analogous conjecture in our context would be that for a cycle $C(t)$ of $f$ the Abelian integral $\int_{C(t)}g$ vanishes if and only if there exist decompositions $f=f_0\circ h$ and $g=g_0\circ h$ such that the projected cycle $h(C(t))$ is trivial. A simple counter-example was given by Pakovich~\cite{P}. 

This allows the formulation of a weaker composition conjecture: The integral $\int_{C(t)}g$ vanishes if and only if there exist decompositions $f=f_k\circ h_k$, $g=g_1\circ h_1+\cdots+g_s\circ h_s$ such that the projected cycles $h_i(C(t))$, $i=1,\ldots,s$ are trivial. In \cite{PM} the authors solved the moment problem in the polynomial case. Their result can be interpreted as saying that the weak composition conjecture is valid for special cycles appearing in the polynomial moment problem. See Section \ref{sec:moment} for more details. We show in Example~\ref{ex:x4} that even the weak composition conjecture is not true for general cycles in the polynomial case.

\section{Main Results}\label{sec:main}

We give a complete solution of  Problem~\ref{problem1} in Theorem~\ref{theo:1} and a partial
solution of Problem~\ref{problem2}  in a generic polynomial case in Theorem~\ref{theo:2}
and in some exceptional cases in Section \ref{balanced}.
Our first result is the following:

\begin{theo}\label{theo:1}
Given $f,g\in\mathbb{C}(z)$, the following conditions are equivalent:
\begin{enumerate}
\item $\int_{C(t)}g\equiv 0$ for every (simple) cycle $C(t)$.
\item There exists $g_0\in \mathbb{C}(z)$,  such that $g=g_0\circ f$.
\end{enumerate}
\end{theo}

\medskip
In order to formulate the main theorem, Theorem~\ref{theo:2}, we
must define the notions of unbalanced, totally unbalanced and
projected cycles. Let $f$ be a polynomial of degree $m$ and
$\Sigma\subset \C$ its set of critical values. Then $f:\C\setminus
f^{-1}(\Sigma)\to\C\setminus\Sigma$ defines a fibration, with
fiber $f^{-1}(t)$ consisting of $m=\deg(f)$ points. Take a base
point $t_0\in\C\setminus\Sigma$. Then this fibration defines a
mapping from the first homotopy group
$\pi_1(\C\setminus\Sigma,t_0)$ to the group of automorphisms
$Aut(f^{-1}(t_0))$. Its image is called the \emph{monodromy group} $G_f$.
It acts transitively on the fiber $f^{-1}(t_0)$ (see, e.g., \cite{Z}).
Changing the base point $t_0$ conjugates all elements of the
group, so the abstract group $G_f$ does not depend on the choice
of the base point $t_0$. Consider the loop in
$\pi_1(\C\setminus\Sigma,t_0)$ winding anti-clockwise around all
critical values in $\Sigma$. It corresponds to a permutation cycle
$\tau_\infty$ of $G_f$ of order $m$. We label the roots
$z_i(t_0)$, $i=1,\ldots,m$ of $f(z_0)=t$ so that this permutation
shifts the indices of the roots by one. The choice of the first
root is arbitrary.

Let $\Gamma_m(f)\subset G_f$ denote the conjugacy class of $\tau_\infty$, that is,
the set of all $\sigma \circ \tau_\infty \circ \sigma^{-1}$, for any $\sigma\in G_f$. Note that $\Gamma_m(f)$ is the set of permutations induced by all paths winding counter-clockwise once around infinity.

\begin{defi}\label{defi:balanced}
We say that a cycle
$C(t)$ of $f$ is \emph{balanced} if
\[
\sum_{i=1}^m n_{p_i}\epsilon_m^i=0,\quad \text{for every }\sigma=(p_1,p_2,\ldots,p_m)\in \Gamma_m(f),\quad 
\]
where $\epsilon_m$ is any primitive $m$-th root of unity.
If $C(t)$ is not balanced, we say that $C(t)$ is \emph{unbalanced}.

\medskip
For each $f_0,h\in \mathbb{C}[z]$ such that $f=
f_0\circ h$ and each cycle $C(t)$ of $f$ we define a cycle
${h}({C}(t))$ of ${f}_0$ called the \emph{projected cycle}
by
\begin{equation}\label{def:cicloproyectado}
{h}(C(t))=\sum_{{h}(z_i(t)) }
\left(\sum_{{h}(z_j)={h}(z_i)}
n_j\right){h}(z_{i}(t))=\sum_{j=1}^d \tilde{n}_j w_j(t).
\end{equation}
Here $w_1(t),\ldots w_d(t)$ are all the different roots
${h}(z_i(t))$ of $f_0(z) = t$.

We say that $C(t)$ is \emph{totally unbalanced} when for every
$f_0,h$ such that $f=f_0\circ h$, the
projected cycle $h({C})$ is unbalanced or trivial.
\end{defi}

The notion of balanced cycle
$C(t)\in \tilde{H}_0(f^{-1}(t))$ is 
related to the behaviour
of the branches $z_i(t)$ when $t$ tends to infinity.
It is well defined, i.e., independent on the way how permutation cycles are written. Indeed,
if we write $\sigma=(p_2, \ldots, p_m, p_1)$, then
\[
\left(n_{p_2}\epsilon_m+n_{p_3}\epsilon_m^2+\ldots+n_{p_m}\epsilon_m^{m-1}+n_{p_1}\right)  \epsilon_m =
\sum_{i=1}^m n_{p_i}\epsilon_m^i.
\]
Thus, $\sum_{i=1}^m n_{p_i}\epsilon_m^i=0$ if and only if
\[n_{p_2}\epsilon_m+n_{p_3}\epsilon_m^2+\ldots+n_{p_m}\epsilon_m^{m-1}+n_{p_1}=0.\]

Note that in particular totally unbalanced cycles are unbalanced as seen from the trivial
decomposition $f=f\circ\id$.

\begin{theo}\label{theo:2}
Let $f\in\mathbb{C}[z]$.

\begin{enumerate}
\item If $C(t)$ is a totally unbalanced cycle, then
\begin{equation}\label{eq:nullAbel}
\int_{C(t)} g\equiv 0
\end{equation}
if and only if there exist $f_1,\ldots,f_s,g_1,\ldots,g_s,h_1,\ldots,h_s\in\mathbb{C}[z]$
such that $f=f_k\circ h_k$, $g=g_1\circ h_1+\ldots+g_s\circ h_s$ and  the projected cycles $h_k(C)$ are trivial
for every $k=1,2,\ldots,s$.

\item If $C(t)$ is an unbalanced cycle, then $g\in\mathbb{C}[z]$
is a solution of \eqref{eq:nullAbel} if and only if there exist
$f_1,\ldots,f_s,g_1,\ldots,g_s,h_1,\ldots,h_s\in\mathbb{C}[z]$
such that $f=f_k\circ h_k$, $g=g_1\circ h_1+\ldots+g_s\circ h_s$,
and for every $k$ 
\[
\int_{{h_k}({C}(t))} g_k \equiv 0,
\]
and the projected cycle ${h_k}({C})$ is
trivial or  balanced.
\end{enumerate}
\end{theo}
In Example \ref{ex:x4}, we show that claim (1) in the Theorem is not true for general cycles contrary to the situation in polynomial moment problem. In order to completely solve Problem~\ref{problem2} for Abelian integrals, it remains to solve it for balanced cycles. Balanced cycles are exceptional cycles of codimension at least one. Moreover, in Section~\ref{balanced} we show that if $m$ is prime, then there are no balanced cycles. We also completely solve the problem for $f(z)=z^m$.

\medskip
In addition to the  problem of persistence of centers, another motivation for our study is the relation of
Problem~\ref{problem2} to the moment problem. The polynomial moment problem consists for a given polynomial $f$ in searching for all polynomials $g$ such that
\[
\int_0^1 f^k(z)g'(z)\,dz=0,\quad \text{for every }k\in\mathbb{N}.
\]
In Section~\ref{sec:moment} we show that the moment problem is equivalent to solving Problem~\ref{problem2} \emph{for a concrete cycle} $C_f(t)$. We show that this cycle is always totally unbalanced, so the first part of Theorem~\ref{theo:2} extends Pakovich and Muzychuk's theorem \cite{PM} solving the polynomial moment problem. In \cite{PRitt} Pakovich gives an explicit solution of the polynomial moment problem. In particular he shows that the function $g$ can be written as a sum of at most three reducible solutions (a polynomial $g$ is a \emph{reducible solution} if there exists $f_0,g_0,h$ such that $\int_{h(C(t))}g_0\equiv 0$, $f=f_0\circ h$ and $g=g_0\circ h$). In Example~\ref{cicloTrivializante} we show that the analogous statement for the tangential center problem is not even true for general totally unbalanced cycles.

\section{Displacement function along a family of zero-cycles}\label{sec:3}

Zero-dimensional Abelian integrals $\int_{C(t)}g$, $C(t)\in \tilde H_0(f^{-1}(t))$ are defined for a couple $f$ and $g$ of polynomials in one variable. 

The definition is analogous to the definition of one-dimensional Abelian integrals $\int_{\gamma(t)}G$, $\gamma(t)\in H_1(F^{-1}(t)$, defined for a polynomial $F$ and a polynomial differential $1$-form $\eta=G_1dx+G_2dy$ in two variables. In this section we push further the analogy between systems depending on one variable and planar systems.  Assume that some point $p\in\R^2$ is a strict local extremum of a real polynomial function $F[x,y]$. Then the system $dF(x,y)=0$ has a center at $p$, parametrized by a family of one-cycles $\gamma(t)\in F^{-1}(t)$.
Consider the deformation
\begin{equation}\label{deformation1}
dF+\epsilon \eta=0.
\end{equation}
For most deformations $\eta$, the center will be broken. A natural question is for what deformations $\eta$, the center is preserved? This is the infinitesimal center problem.

The problem is tackled by considering the displacement function on a transversal to the family of cycles 
$\gamma(t)$ parametrized by $F$. The displacement function $D_\epsilon(t)$ measures the displacement along a trajectory of the deformed system \eqref{deformation1} from the starting point on the transversal to the first return point to the transversal. A center is preserved if $D_\epsilon\equiv0$. A classical result is that Abelian integrals give the first order term of the displacement function:
\begin{equation}\label{D}
D_\epsilon(t)=-\epsilon\int_{\gamma(t)}\eta+o(\epsilon).
\end{equation}
Of course the vanishing of the first order (Abelian integral) term is a necessary, but not sufficient condition for the persistence of a center. It is named the tangential (or first order) center problem.

The above problem can be complexified. One considers the foliation defined by \eqref{deformation1}.
Taking a transversal and lifting each loop $\gamma(t_0)\subset F^{-1}(t)$ to nearby leaves the holonomy associated to $\gamma(t_0)$ is defined. By assumption, for $\epsilon=0$ all holonomies are identity (i.e., the associated displacement functions are zero-functions). 

\medskip

Consider now the problem in one-dimensional space. Let $f\in\C(z)$ be a function and $C(t)\in\tilde H_0(f^{-1}(t))$ a family of cycles of $f$. Consider a perturbation of $f$ by $g$
$$
f+\epsilon g
$$
and the deformed cycle $C_{\epsilon}(t)=\sum_{i=1}^m n_i z_i(t,\epsilon)$, where $(f+\epsilon g)(z_i(t,\epsilon))=t$, for $t \in \C \backslash \Sigma_\epsilon$, where $\Sigma_\epsilon$ is the set of critical values of $f+\epsilon g$. We define the \emph{displacement function of $f+\epsilon g$ along the cycle $C(t)$} by
\[
\Delta_\epsilon(t)=\sum_{i=1}^m n_i f(z_i(t,\epsilon))=\int_{C_\epsilon(t)}f.
\]
As in the case of one-dimensional loops, zero-dimensional Abelian integrals correspond to the principal part of the displacement function of the perturbation $f+\epsilon g$  along a family of zero-dimensional cycles. Indeed, from $(f+\epsilon g)(z_i(t,\epsilon))=t$, it follows that $\sum_{i=1}^m n_i (f+\epsilon g)(z_i(t,\epsilon))=0$, giving
\begin{equation}\label{eq:displacement}
\Delta_\epsilon(t)=-\epsilon\int_{C_\epsilon(t)} g=-\epsilon\int_{C(t)} g+o(\epsilon).
\end{equation}

Problems \ref{problem1} and \ref{problem2} can be seen as first order (called also tangential) problems of analogous problems  for the displacement function along  families of zero-cycles:

\begin{problem}\label{problem1'}
Characterize rational functions $f$, $g$  such that for any family of cycles $C(t)$ of $f$ the displacement function $\Delta_\epsilon$ of $f+\epsilon g$, along any zero cycle $C(t)$ vanishes identically.
\end{problem}

\begin{problem}\label{problem2'}(Infinitesimal center problem)
Characterize rational functions $f$, $g$ and cycles $C(t)$ of $f$
such that the displacement function $\Delta_\epsilon$ of $f+\epsilon g$, along the family of zero-cycles $C(t)$ vanishes identically.
\end{problem}

Solution of Problem \ref{problem2'} is a corollary to Theorem \ref{theo:1}.

\begin{coro}\label{coro:1}
Given $f,g\in\mathbb{C}(z)$, the following conditions are equivalent:
\begin{enumerate}
\item $\Delta_\epsilon(t)\equiv 0$ for every (simple) cycle $C(t)$.
\item There exists $g_0\in \mathbb{C}(z)$,  such that $g=g_0\circ f$.
\end{enumerate}
\end{coro}

Problem \ref{problem2'} would probably require studying vanishing of some {\it iterated integrals} on zero-cycles. We hope to address this problem in the future.

\section{Imprimitivity systems of the monodromy group and equivalence classes of decompositions}\label{sec:preliminaries}

The aim of this section is to recall some definitions and results
about the monodromy group and to study the relationship between the
decompositions of a rational function and its monodromy group.

Given a rational function $f\in\mathbb{C}(z)$ of degree $m$, let
$\Sigma\subset \mathbb{C}$ be the set of critical values of
$f$. The monodromy group $G_f$ is defined as in the polynomial case and acts transitively on generic fibers $f^{-1}(t),$ $t\in\mathbb{C}\setminus\Sigma$.

As shown in \cite{F} the monodromy group $G_f$ is the Galois group of the Galois extension of
$\mathbb{C}(t)$ by the $m$ preimages $z_1(t), \ldots, z_m(t)$ of $t$ by $f$, that is,
\[
G_f = Aut(\mathbb{C}(z_1, \ldots, z_m)/\mathbb{C}(t)).
\]

Let $X = \{ 1, \ldots, m\}$ and let $G \subseteq S_m$ be a
transitive permutation group on $X$. A subset $B \subseteq X$ is
called a {\em block} (\cite{Wielandt}) of $G$ if for each $g \in
G$ the image set $g(B)$ and $B$ are either equal or disjoint.
Given a block $B$, the set $\mathcal{B} = \{ g(B) | g \in G\}$
forms a partition of $X$ into  disjoint blocks of the same
cardinality which is called an {\em imprimitivity system} (or a
{\em complete block system}). Each permutation group $G \subseteq
S_m$ has two {\em trivial} imprimitivity systems: $\{X\}$ and
$\{ \{x\} \}_{x \in X} $.

Given $f(z) \in \mathbb{C}(z)$ we say that two {\em
decompositions} of $f$, $f=f_1 \circ h_1$ and  $f=f_2 \circ h_2$,
are {\em equivalent} if $h_1$ and $h_2$ define the same field
over every preimage of $t$ by $f(z)$, that is, if
$\mathbb{C}(h_1(z_i))=\mathbb{C}(h_2(z_i))$ for every $i=1,\ldots,m$.

\begin{prop}\label{prop:blocks}
Let $f(z) \in \mathbb{C}(z)$. There exists a one-to-one
correspondence between imprimitivity systems of $G_f$ and
equivalence classes of decompositions of $f$.

Moreover, if $f=f_0 \circ h$, $\mathcal{B}_h$ is the
imprimitivity system corresponding to $h$ and $i\in B\in\mathcal{B}_h$, then
$B=\{k\colon h(z_k)=h(z_i)\}$.
\end{prop}

\begin{rema}\label{remark:blocks}
Let us point out that the second statement in the proposition is
in particular telling us that given a cycle $C(t)= \sum_{i=1}^m
n_i z_i(t)$ of $f$, the projected cycle $h(C(t))$ of $f_0$
given by the expression (\ref{def:cicloproyectado}) can be written
as
\[
h(C(t))= \sum_{h(z_i(t))} \left( \sum_{j \in B_i} n_j \right) h(z_i(t)),
\]
where $B_i$ is the block containing $i$ of the imprimitivity
system $\mathcal{B}$ corresponding to the decomposition $f=\tilde
f \circ h$.
\end{rema}

\begin{proof}
Let $f(z)$ be a rational function of degree $m$, and let us denote
by $z_1(t), \ldots, z_m(t)$ the preimages of $t\in \mathbb{C} \backslash
\Sigma$ as above.

Let $f = f_0 \circ h$ be a decomposition and let us find the
associated imprimitivity system $\mathcal{B}_h$. To do so, we only
have to find a block $B$, since by definition, $\mathcal{B}_h
 := \{ \sigma(B) : \sigma \in G_f \}$ is an imprimitivity system of $G_f$. For any preimage $z_i$ of $t \in
\mathbb{C} \backslash \Sigma$, we have that
\[
t = f(z_i) = f_0 (h(z_i)),
\]
and therefore
\[
\mathbb{C}(t) \subseteq \mathbb{C}(h(z_i)) \subseteq \mathbb{C}(z_i) \subseteq \mathbb{C}(z_1, \ldots, z_m).
\]
By the Fundamental Theorem of Galois theory (in Artin's version,
see \cite{Artin, W}) there exists a subgroup $H$ of the monodromy
group $G_f$ containing $H_i$, the stabilizer of $i$, such that
\[
\mathbb{C}(h(z_i)) = \mathbb{C}(z_1, \ldots, z_m)^H.
\]
As a consequence, $h(z_i) = \tau(h(z_i)) = h(\tau(z_i))$ for every
$\tau \in H$. The orbit $B$ of $i$ by the action of this subgroup
$H$ is the block we are looking for:
\[
B := \{ \sigma(i) : \sigma \in H \} \simeq H/H_i.
\]
It is obvious that $\tau(B)=B$ for every $\tau \in H$, and it is
also immediate to see that $\sigma(B) \cap B = \emptyset$ when
$\sigma \in G_f \backslash H$. Then, we also have proved that
$h(z_i)$ is constant on $i \in B \in \mathcal{B}_h$. Conversely,
if $h(z_i)=h(z_k)$, then $k \in B$: since the monodromy group $G_f$ is
transitive, then there exists some $\sigma \in G_f$ such that
$\sigma(i)=k$, and therefore
\[
\sigma(h(z_i)) = h(\sigma(z_i)) = h(z_k) = h(z_i),
\]
which implies that $h(z_i)$ is invariant by $\sigma$. Then
$\sigma$ must belong to $H$ and $k$ belongs to $B$. For another
different block $\tilde B = \sigma(B)$, a similar argument proves the
same statement. Then, the blocks of the imprimitivity system
$\mathcal{B}_h$ verify that $h(z_i)=h(z_j)$ if and only if $i,j$
belongs to the same block of $\mathcal{B}_h$, proving that the
imprimitivity system we get does not depend on the choice of
$z_i$.

Now let $\mathcal{B}$ be an imprimitivity system for $G_f$. The
block $B$ containing the element $i$ provides a subgroup of $G_f$
containing the stabilizer $H_i$ of $i$. Precisely,
\[
G_B := \{ \tau \in G_f : \tau(B)= B \}.
\]
By the Fundamental Theorem, the following inclusions of groups
\[
\{ {\rm Id} \} \subseteq H_i \subseteq G_B \subseteq G_f
\]
yield the following inclusions of fields
\[
\mathbb{C}(t) \subseteq L_B: = \mathbb{C}(z_1, \ldots, z_m)^{G_B}
\subseteq \mathbb{C}(z_i) \subseteq \mathbb{C}(z_1, \ldots, z_m).
\]
By L\"{u}roth's Theorem \cite{W} applied to $\mathbb{C} \subset
L_B \subseteq \mathbb{C}(z_i)$, $L_B = \mathbb{C}(h(z_i))$ for
some rational function $h$. Since $f(z_i)=t \in L_B$, there exists
another rational function $f_0$ such that $f = f_0 \circ
h$. Moreover, $h(z_i)=h(z_j)$ for any $i,j\in B$.

If we choose an element $k$ of a different block $\tilde B$, then we get
another subgroup $G_{\tilde B}$ of $G_f$, which contains the stabilizer
$H_k$ of $k$. Again by the Fundamental Theorem, it provides a
field $L_{\tilde B} := \mathbb{C}(z_1, \ldots, z_m)^{G_{\tilde B}}$, which by
L\"{u}roth's Theorem is generated by some rational function $\tilde h$,
that is, $L_{\tilde B} = \mathbb{C}(\tilde h(z_k))$. Therefore, there exists
another rational function $\tilde f_0$ such that $ f= \tilde f_0
\circ \tilde h$. Let us see how different these functions $h$ and $\tilde h$
are. Since $G_f$ is transitive, there exists a permutation $\tau$
such that $\tau(i)=k$. Therefore, $\tau(B)=\tilde B$ (or
$B=\tau^{-1}(\tilde B)$), and it is immediate to check that
\begin{equation}\label{grupos_conjugados}
\tau \circ G_B \circ \tau^{-1} = G_{\tilde B}.
\end{equation}
Since by the Fundamental Theorem we know that $G_B$ and
$G_{\tilde B}$ are the Galois groups of $\mathbb{C}(z_1,\ldots,z_m)$
over $\mathbb{C}(h(z_i))$ and $\mathbb{C}(\tilde h(z_k))$ respectively,
let us see what equality \eqref{grupos_conjugados} gives: an
element $\sigma' \in G_{\tilde B}$ is an algebra automorphism of
$\mathbb{C}(z_1,\ldots,z_m)$ such that $\sigma'(\tilde h(z_k))=\tilde h(z_k)$.
Since $\sigma' = \tau \circ \sigma \circ \tau^{-1}$ for some
$\sigma \in G_B$, that means
\[
\begin{matrix}
\tau ( \sigma ( \tau^{-1}(\tilde h(z_k)))) & = & \tilde h(z_k) \\
\parallel & & \parallel \\
\tau ( \sigma ( \tilde h(z_i))) &  & \tau ( \tilde h(z_i)),
\end{matrix}
\]
which implies that $\sigma ( \tilde h(z_i)) = \tilde h(z_i)$. Also, as $\sigma
= \tau^{-1} \circ \sigma' \circ \tau $ leaves invariant $h(z_i)$,
it holds that
\[
\begin{matrix}
\tau^{-1} ( \sigma' ( \tau(h(z_i)))) & = & h(z_i) \\
\parallel & & \parallel \\
\tau^{-1} ( \sigma' ( h(z_k))) &  & \tau^{-1} ( h(z_k)),
\end{matrix}
\]
which implies that $\sigma'(h(z_k)) = h(z_k)$. To sum up, equality
\eqref{grupos_conjugados} implies that $G_B$ is not only the
Galois group of $\mathbb{C}(z_1,\ldots,z_m)$ over
$\mathbb{C}(h(z_i))$, but also over $\mathbb{C}(\tilde h(z_i))$.
Similarly, $G_{\tilde B}$ is not only the Galois group of
$\mathbb{C}(z_1,\ldots,z_m)$ over $\mathbb{C}(\tilde h(z_k))$, but also
over $\mathbb{C}(h(z_k))$. By the Fundamental Theorem, we
obtain that $L_B= \mathbb{C}(h(z_i)) = \mathbb{C}(\tilde h(z_i))$ and
$L_{\tilde B} = \mathbb{C}(\tilde h(z_k)) = \mathbb{C}(h(z_k))$. Therefore,
given an imprimitivity system $\mathcal{B}$, we get not only one
decomposition of $f$, but a whole equivalence class of
decompositions.
\end{proof}

\section{Zero-dimensional Abelian integrals and displacement functions vanishing on any cycle}\label{sec:all}

The aim of this section is to solve Problems \ref{problem1} and \ref{problem1'}  proving Theorem~\ref{theo:1} and Corollary \ref{coro:1}. First we give a simple sufficient condition for the vanishing of the displacement function $\Delta_\epsilon$ (and hence also the corresponding Abelian integral) of the family $f+\epsilon g$, $f,g\in \mathbb{C}(z)$, along a family of cycles $C(t)$ of $f$.

\begin{prop}\label{prop:compositionconjecture}
Assume that $f=f_0\circ h$, $g=g_0\circ h$, where $f_0,g_0,h\in \mathbb{C}(z)$.
Let $C(t)=\sum_k n_k z_k(t)$ be a cycle of $f$ and $C_\epsilon(t)=\sum_k n_k z_k(t,\epsilon)$
its continuation as a cycle of $f+\epsilon g$. Let
\[
h(C(t))=\sum_{h(z_k(t))} \left( \sum_{h(z_j(t))=h(z_k(t))}n_j \right) h(z_k(t))
\]
and
\[
h(C_\epsilon(t)) = \sum_{h(z_k(t,\epsilon))} \left(
\sum_{h(z_j(t,\epsilon))=h(z_k(t,\epsilon))}n_j \right)
h(z_k(t,\epsilon))
\]
be the cycle of $f_0$ and $f_0+\epsilon g_0$ obtained by projection of $C(t)$ and $C_\epsilon(t)$ by $h$ respectively. Then,
\begin{enumerate}
\item $\Delta_\epsilon(t)=\int_{h(C_\epsilon(t))}f_0$ and
$\int_{C(t)}g=\int_{h(C(t))}g_0$.
\item
In particular, if the projected cycle $h(C(t))$ is trivial, then $h(C_{\epsilon}(t))$ is trivial,
\[
\Delta_\epsilon(t)\equiv 0
\]
and
\[
\int_{C(t)} g\equiv0.
\]
\end{enumerate}
\end{prop}

\begin{rema}
In other words, the assumption (2) of the proposition can be reformulated by saying that the cycle $C(t)$ projects by $h$ to a trivial cycle
\[
h(C(t))=\sum_{f_0(w_i(t))=t} \left(\sum_{h(z_k(t))=w_i(t)} n_k \right)  w_i(t)
\]
of $f_0$ and, moreover, the function $g$ factors through the same $h$, so a change of coordinates leading to integration on the trivial cycle is possible.
\end{rema}

\begin{proof}(1) Let $C(t)=\sum_k n_k z_k(t)$ be a cycle and assume that $f=f_0\circ h$, $g=g_0\circ h$.
Set $w_i(t,\epsilon)=h(z_i(t,\epsilon))$. Then
\[
\begin{split}
\Delta_\epsilon(t) & = \sum_{i=1}^m n_i f(z_i(t,\epsilon))=\sum_{h(z_i(t,\epsilon))}
\left(\sum_{h(z_j(t,\epsilon))=h(z_i(t,\epsilon))} n_i\right)f_0(w_j(t,\epsilon)) \\ & =\int_{h(C_\epsilon(t))}f_0.
\end{split}
\]
The second claim of (1) follows the same way.

(2) If $h(C(t))$ is a trivial cycle, then all its coefficients are equal to zero. The dependence with respect to $\epsilon$ being continuous, it follows that the cycle $h(C_\epsilon(t))$ is trivial too. Now the claim follows from (1).
\end{proof}

Next result shows that in case of simple cycles condition (2) of Proposition \ref{prop:compositionconjecture} is necessary and sufficient for the vanishing of Abelian integrals.

\begin{prop}\label{prop:simplecycle}
Let $f(z) \in \mathbb{C}(z)$, and consider a simple cycle $C(t)=z_i(t)-z_j(t)$. Then
\[
\int_{C(t)}g\equiv 0
\]
if and only if there exist $f_0,g_0,h\in\mathbb{C}(z)$ such that  $f=f_0\circ  h$, $g=g_0\circ h$ and $ h(z_i(t)) = h(z_j(t))$.
\end{prop}

\begin{rema}\label{rema:simplecicles}
Let us observe that since $h(z_i(t))= h(z_j(t))$, by Proposition \ref{prop:blocks} the set $\{i,j\}$ is included in some block $B$ of the imprimitivity system $\mathcal{B}_{ h}$ corresponding to the decomposition
$f=f_0\circ h$.

Proposition~\ref{prop:simplecycle} essentially follows from Theorem~4 of \cite{GM} and L\"uroth's Theorem, but we include the proof for completeness.
\end{rema}

\begin{proof}
The sufficient condition is a consequence of Proposition
\ref{prop:compositionconjecture}.

For the necessary condition, let us define
\[
\mathcal{W}=\{h\in\mathbb{C}(z)\colon h(z_i(t))=h(z_j(t))\}.
\]
Since
\[
\mathbb{C}\subset \mathcal{W}\subset \mathbb{C}(z),
\]
by L\"{u}roth's Theorem $\mathcal{W}=\mathbb{C}(h(z))$ for some rational function $h$. Since $f,g\in\mathcal{W}$, then $f=f_0\circ h$, $g=g_0\circ  h$ for some $f_0, g_0 \in \mathbb{C}(z)$.
\end{proof}

Now, we are ready to prove Theorem~\ref{theo:1}.

\begin{proof}[Proof of Theorem~\ref{theo:1} and Corollary \ref{coro:1}]
By Proposition~\ref{prop:compositionconjecture}, in Theorem~\ref{theo:1} $(2)$ implies $(1)$, similarly in Corollary \ref{coro:1} $(2)$ implies $(1)$. Therefore, to conclude it is sufficient to prove that  in Theorem~\ref{theo:1} $(1)$ implies $(2)$.

Let $g\in\mathbb{C}(z)$ be such that $\int_{C(t)}g\equiv 0$ for every (simple) cycle $C(t)$. First, let us observe that the condition $\int_{C(t)} g \equiv 0$ for every (simple) cycle $C(t)$ is equivalent to
\[
g(z_1(t))=g(z_2(t))=\ldots=g(z_m(t)).
\]
Therefore, $g(z_i) \in \mathbb{C}(z_1, \ldots, z_m)$ is invariant under the action of the whole monodromy group $G_f$, that is, $g(z_i) \in \mathbb{C}(z_1, \ldots, z_m)^{G_f} = \mathbb{C}(t)$, which means that
there exists $g_0(z) \in \mathbb{C}(z)$ such that $g(z_i(t)) = g_0(t) = g_0(f(z_i(t)))$. Then $g = g_0 \circ f$.
\end{proof}

Under a generic condition on $f$ we show that the existence of a persistent center or a tangential center for any non-trivial cycle $C(t)$ is equivalent to the integrability of the system as characterized in Theorem~\ref{theo:1}. Recall that a point $z_0$ is of Morse type if $f(z)=f(z_0)+a(z-z_0)^2+\circ(z-z_0)^2$,
$a\in\mathbb{C}$. Generically all critical points are of Morse type and all critical values are different.

\begin{prop}
Given $f \in \mathbb{C}(z)$, assume that every critical value corresponds to only one critical point, and moreover that they are all of Morse type. For any non-trivial cycle $C(t)$ the following conditions are equivalent:
\begin{enumerate}
 \item $\int_{C(t)}g\equiv 0$ .
 \item $\Delta_\epsilon(t)\equiv 0$ along $C(t)$.
 \item There exists $g_0\in \mathbb{C}(z)$ such that $g=g_0\circ f$.
\end{enumerate}
\end{prop}

\begin{proof}
We prove that $(1)$ implies $(3)$ for a given cycle $C(t)$, since, by Theorem~\ref{theo:1}, $(3)$ implies $(1)$ and $(2)$, and trivially $(2)$ implies $(1)$.

The permutation in the monodromy group associated with a critical value with a Morse point is a transposition.
Let the cycle $C(t)$ be given by \eqref{cycle}.

First, we shall prove that there is a transposition $\tau_s$ corresponding to a critical value $s$ and $k\in\{1,2,\ldots,m\}$ such that $n_k\neq n_{\tau_s(k)}$. Note that otherwise, $n_k= n_{\tau_s(k)}$ for every $k\in\{1,2,\ldots,m\}$ and every critical value $s$. Since the monodromy group is transitive and it is generated by the permutations corresponding to critical values, we would have $n_k=a$ for every $k\in\{1,2,\ldots,m\}$, and certain $a\in \mathbb{Z}$. Now, since $C(t)$ is a cycle, $0=\sum n_k=m a$ and the cycle $C(t)$ would be trivial.

Take $s$ and $k\in\{1,2,\ldots,m\}$ such that $n_k\neq n_{\tau_s(k)}$. By analytic continuation, 
\[
\sum_{i=1}^m n_{\tau_s(i)} g(z_i(t))\equiv 0.
\]
Therefore,
\[
0 = \sum_{i=1}^m n_i g(z_i(t))-\sum_{i=1}^m n_{\tau_s(i)}
g(z_i(t))=\left(n_k-n_{\tau_s(k)}
\right)\left(g(z_k(t))-g(z_{\tau_s(k)}(t))\right),
 \]
so the integral $\int_{z_k(t)-z_{\tau_s(k)}(t)}g$ vanishes. By Proposition~\ref{prop:simplecycle} and Remark~\ref{rema:simplecicles}, there exists $f_0, g_0, h$ such that $\{k,\tau_s(k)\}\subset B\in\mathcal{B}_{h}$, $f=f_0\circ h$ and $g=g_0\circ h$.

If $\tau=(i,j)\in G_f$ and $\mathcal{B}$ is a non-trivial imprimitivity system of $G_f$, then $\{i,j\}\subset B$ for some block $B\in \mathcal{B}$. Since $G_f$ is generated by the transpositions corresponding to the critical values and it is transitive, then for each $\{i,j\}\subset \{1,\ldots,m\}$, there
is a chain of transpositions $(i,i_1)$, $(i_1,i_2)$, $\ldots$, $(i_r,j)$. Therefore, there is only one block $B=\{1,\ldots,m\}$, and the unique imprimitivity systems of $G_f$ are the trivial ones. Thus, $\mathcal{B}_{h}$ must be the trivial imprimitivity system with a unique block and, in consequence, $h=f$.
\end{proof}

If the polynomial $f$ is not generic, then there are some solutions $g$ which are not of the form $g=g_0\circ f$. Moreover, there exists solutions that are not of the form $g=g_0\circ h$ with $f=f_0\circ h$.

\begin{exam}
Let $f(x)=x^4-x^2$. Let $x_1(t)<x_2(t)<x_3(t)<x_4(t)$,
$t\in(-1/4,0)$, be the real roots of $f$. Consider the cycle
\[C(t)=x_1(t)-x_2(t)-x_3(t)+x_4(t).
\]
Here $t$ can be analytically extended to $\C\setminus\{-1/4,0\}$.
Then for any odd function $g$ it follows $\int_{C(t)}g\equiv0$,
because $f(x)=t$ is a biquadratic equation and therefore
$x_1(t)=-x_4(t)$, $x_2(t)=-x_3(t)$. Take, for instance, $g(x)=x^3$.

Let us show that there not exist $f_0,g_0,h$ with $deg(h)>1$ such that $f=f_0\circ h$, $g=g_0\circ h$.
Indeed, assume that there exist such decompositions and let $\mathcal{B}_h$
be the associated imprimitivity system, and $B\in \mathcal{B}_h$. Since the
imprimitivity system is not trivial, then $B=\{i,j\}$ and $h(x_i)=h(x_j)$. On the other hand, $g$ does not
admit decompositions, so $g=h$, but $g(x_i)\neq g(x_j)$ for every $i\neq j$.
\end{exam}

\section{Tangential center problem for polynomials}\label{sec:t2}

In this section we focus on the case when $f,g$ are polynomials.
We shall prove Theorem~\ref{theo:2} and some other results under
the assumption that the cycle $C(t)$ given by \eqref{cycle} is
totally unbalanced.

Consider a cycle $C(t)$ given by \eqref{cycle}. It is represented by the vector $\bar n=(n_1,\ldots,n_m)$.
Let $\sigma\in G_f$ be an element of the monodromy group of $f$.
It transforms the cycle $C(t)$ into the cycle
\begin{equation}\label{eq:accion}
\sigma(C(t))=\sum_{i=1}^m n_{i}  z_{\sigma(i)}(t)=\sum_{i=1}^m n_{\sigma^{-1}(i)}z_i(t).
\end{equation}

We consider the action of
$G_f$ on $\mathbb{C}^m$ defined by
\[
\sigma(v) = (v_{\sigma(1)},\ldots,v_{\sigma(m)}),\quad \text{for
every }v=(v_1,\ldots,v_m)\in \mathbb{C}^m,\ \sigma\in G_f.
\]
Replacing $\sigma^{-1}$ by $\sigma$, we see that the orbit of the cycle $C(t)$ by the action of the monodromy group $G_f$ is represented by the space $V\subset\mathbb{C}^m$ defined by
\[
V=<\sigma(\bar n)\colon \sigma\in G_f>.
\]

Let us choose $\tau\in \Gamma_m(f)$ and
number the branches in such a way that $\tau=(1,2,\ldots,m)$.

Let $\mathcal{B}$ be an imprimitivity system of $G_f$. By
Proposition \ref{prop:blocks} it corresponds to a class of
decompositions $f=f_0\circ h$, where $f_0,h$ can be chosen
polynomials (see Lemma~3.5 of \cite{CM}). Let $d=\deg h$. Then
$d|m$ and $m/d=\deg f_0$. Since $(1,2,\ldots,m)\in G_f$, then each
$B_i\in\mathcal{B}$, $i=1,\ldots,m/d$ consists of the congruence
class of $i$ modulo $d$ in $\{ 1,2,\ldots, m\}$. For each
decomposition as above, the
corresponding projected cycle is of the form 
\[
\tilde C_d(t)=\sum_{1\leq i\leq m/d } \left(\sum_{j \equiv_{d} i }
n_j\right) w_{i}(t),
\]
where $w_i(t)$ are all the different values of $h(z_j(t))$.

Let $D(f)$ denote the set of divisors $r$ of $m=\deg f$
such that there is a decomposition $f=f_k\circ h_k$
with $\deg(h_k)=r$. Note that in $D(f)$ we have the partial
order induced by divisibility. We say that $r_1,\ldots,r_l$ is a
is a {\em complete set} of elements of $D(f)$ covered by $r$
if $r_k|r$ and $r_1,r_2,\ldots,r_l$ are maximal among divisors of $r$
in $D(f)\backslash\{r\}$.

For every $r \in D(f)$, let us denote by $V_r$ the set of
$r$-periodic vectors.

\begin{lema}[\cite{PM}]\label{lem:irreducible}
Each $G_f$-irreducible subspace of $\mathbb{Q}^m$ has the form
\[
U_r:=V_r\cap\left(V_{r_1}^\bot\cap \ldots \cap
V_{r_l}^\bot\right),
\]
where $r\in D(f)$ and $r_1,\ldots,r_l$ is a complete set of
divisors of $D(f)$ covered by $r$. The subspaces $U_r$ are
mutually orthogonal and every $G_f$-invariant subspace of
$\mathbb{Q}^m$ is a direct sum of some $U_r$ as above.
\end{lema}

\begin{lema}[\cite{PM}]\label{lem:4.3}
Let $g$ be a polynomial and let
\[
g(z_1(t)) = \sum_{
k\geq -\deg(g)}
 s_k t^{-k/m}
\]
denote the Puiseaux expansion at infinity of $g(z_1(t))$, where
$z_1(t)$ is one of the branches of the preimage of the polynomial
$f$. For any $c \in D(f)$, $c \neq m$, if we define
\[
\psi_c(t) = \sum_{
k\equiv_{m/c} 0,
k\geq -\deg(g)}
 s_k t^{-k/m},
\]
there exist $w(z) \in \C[z]$ such that
\[
\psi_c(t)=w(z_1(t)).
\]
Moreover, there exist $f_0, g_0, h \in \C[z]$, with
$deg(h)>1$, such that
\[
f = f_0 \circ h, \quad w = g_0 \circ h.
\]
\end{lema}

Now we are in conditions of proving Theorem~\ref{theo:2}.

\begin{proof}[Proof of Theorem~\ref{theo:2}]

First we prove $(2)$.

Assume that there exist
$f_1,\ldots,f_s,g_1,\ldots,g_s,h_1,\ldots,h_s \in \mathbb{C}[z]$
such that $f=f_k\circ h_k$, $k=1,\ldots,s$, $g=g_1\circ
h_1+\ldots+g_s\circ h_s$. Set $d_k=\deg(h_k)$ and let $h_k(C(t))$
denote the projected cycle
\[
h_k(C(t))=\sum_{1\leq i\leq m/d_{k} } \left(\sum_{j \equiv_{d_{k}}
i } n_j\right) w_{i}(t), \quad k=1,\ldots,s.
\]
Assume that for every $k$, either the projected cycle is trivial
or balanced, and
\[
\int_{h_k(C(t))} g_k \equiv 0.
\]
Then by Proposition~\ref{prop:compositionconjecture}, $g_k\circ h_k$ is a
solution of \eqref{eq:nullAbel} for every $k$. By linearity, $g$ is also a
solution of \eqref{eq:nullAbel}.

\vspace{0.5cm} Now, assume that $g$ is a solution of
\eqref{eq:nullAbel}.
By analytic continuation this means that
\begin{equation}\label{eq:sigma}
\sum_{i=1}^m n_{\sigma(i)}g(z_i(t))\equiv0.
\end{equation}
Assume that $C(t)$ is unbalanced. Thus, there exists a permutation cycle of order $m$, $\tau=(p_1,p_2,\ldots,p_m)\in \Gamma_m(f)$,
such that
\[
\sum_{i=1}^m n_{p_i}\epsilon_m^i\neq 0.
\]
There is no loss of generality in assuming
that $\tau=(1,2,\ldots, m)$. Consider the orbit of the cycle
\[
V=<\sigma(\bar n)\colon \sigma\in G_f>.
\]

Let us define the vectors
\[
w_k=(1,\epsilon_m^k,\epsilon_m^{2k},\ldots,\epsilon_m^{(m-1)k})
\]
for every natural $k$, where $\epsilon_m = e^{2\pi i/m}$. Note
that these vectors for $k=1,\ldots,m$ form an orthogonal basis of
$\mathbb{C}^m$. Moreover, for every $r \in D(f)$, the (complexified) subspace
$V^{\C}_r$ is generated by the vectors $w_k$ such that $(m/r)|k$.

By Lemma~\ref{lem:irreducible}, for $r=m$ we have that
\[
W=U_m=V_{r_1}^\bot\cap \ldots \cap V_{r_l}^\bot
\]
is a $G_f$-irreducible subspace of $\Q^m$, where $r_1,\ldots,r_l$
are the maximal elements of $D(f) \setminus \{m\}$.

Let us prove that $W\subset V$. Since $V,W$ are invariant by the
action of $G_f$ and $W$ is irreducible, then either $W\subset V$
or $W\subset V^\perp$. We show that $w_1\in W^\mathbb{C}$
and $w_1\not\in V^{\mathbb{C}\perp}$, thus $W^\mathbb{C}\not\subset V^{\mathbb{C}\perp}$.
Indeed, for every $i$, $V^{\C}_{r_i}$ is
generated by the vectors $w_k$ where $(m/r_i) | k$. Since $r_i
\neq m$, then $k>1$ and these $w_k$ are orthogonal to $w_1$.
Therefore $w_1 \in W^{\C}$, but since $C(t)$ is unbalanced, $w_1$
is not orthogonal to $\bar{n} \in V$. Then $V$ and $W$ are not
orthogonal, thus showing that $W\subset V$.

Now, take the Puiseux expansions of $g(z_i(t))$ at infinity.
Applying Newton's Polygon Method to $f(z)=t$ (see, e.g., \cite[p.
45]{CH}) $z_i(t)$ can be written as an analytic function in  $\epsilon_m^{i-1}t^{-1/m}$,
and composing with $g$ we get
\[
g(z_i(t))=\sum_{k=-\deg(g)}^\infty s_k \epsilon_m^{(i-1)k}
t^{-k/m}.
\]
Now, \eqref{eq:sigma} yields
\begin{equation}\label{eq:condsk}
s_k \sum_{i=1}^m n_{\sigma(i)} \epsilon_m^{(i-1) k}=0,\quad k\geq
-\deg(g),\ \sigma\in G_f.
\end{equation}

Assume that $s_{k_0}\neq 0$ for some ${k_0} \geq-\deg(g)$. Then
$w_{k_0}$ is orthogonal to $V$, and therefore, to $W$. Since
$W^{\C\perp}$ is generated by the vectors $w_k$, $(m/r)|k$, $r$ maximal
in $D(f)\setminus \{m\}$, this implies that $w_{k_0}$ is a linear
combination of these vectors. But $\{w_i\}_{1, \ldots, m}$ are
linearly independent, therefore $w_{k_0}$ coincides with one of
them. Consequently, $(m/c)|k_0$ for some $c$ maximal element of
$D(f)\backslash \{m\}$.

Let $r=m/c$. Define
\[
\psi_c(t)=\underset{k\geq -\deg(g)}{\underset{k\equiv_r 0}{\sum}} s_k t^{-k/m}.
\]
By Lemma~\ref{lem:4.3}, there exist $f_0,g_0,h,w\in\mathbb{C}[z]$, such that $\psi_c=w\circ z_1$ and 
$f=f_0\circ h$, $w=g_0\circ h$, with $\deg(h)>1$.

The Puiseux expansion of $w(z_i(t))$ at infinity is
\[
w(z_i(t))=\underset{k\geq -\deg(g)}{\underset{k\equiv_r 0}{\sum}} s_k \epsilon_m^{(i-1)k}t^{-k/m}.
\]
Therefore, $w$ satisfies \eqref{eq:condsk} and as a consequence
\[
\int_{C(t)} w = \sum_{i=1}^m n_i w(z_i(t))\equiv 0.
\]
Now, since $w=g_0 \circ h$, projecting the cycle
$C(t)$ by $h$,
\[
\int_{h(C(t))} g_0 = \sum_{1\leq i\leq m/d } \left(\sum_{j \equiv_{d} i } n_j\right)
g_0(w_{i}(t))\equiv 0,
\]
where $w_i(t)=h(z_i)$, and $d$ is the number of different
$w_i(t)$.

Since $g-w$ has $s_{i m/c}=0$ for every $i\geq -\deg(g)$,
by induction on the number of maximal elements $\bar c$ of $D(f)\backslash \{m\}$
such that $s_{i m/\bar c}\neq 0$ for some $i$, we obtain
\[
g=\sum_k g_k,
\]
where the Abelian integral of $g_k$ along the cycle
\[
C_k(t):=\sum_{1\leq i\leq m/d_k } \left(\sum_{j \equiv_{d_k} i } n_j\right) w_{j}(t)
\]
vanishes for every $k$.
It only remains to prove that we can assume that either $C_k(t)$ is trivial or
balanced. If it is none of them, then it is unbalanced, and
applying again the arguments above, $g_k=\sum_j g_{k_j}$, where
now the cycles associated with each $g_{k_j}$ are projections of the
cycle $C_k(t)$. Repeating this argument recursively, one obtains
that either $C_k(t)$ is trivial or balanced.

In particular, if $C(t)$ is totally unbalanced, then every
$C_k(t)$ must be trivial, which proves $(1)$.
\end{proof}

In~\cite{PRitt} Pakovich proves in particular that any solution of the
polynomial moment problem can be expressed as a sum of at most
three reducible solutions. We show that in our case of general
cycles there exist solutions with arbitrarily many terms
which cannot be reduced.

Let $f(z)$ be a polynomial of degree $m$ such that $f^{-1}(t)$ has
$m$ different branches $z_1(t), \ldots, z_m(t)$ for every $t \in
\mathbb{C}$. Let $C(t)$ be a cycle given by (\ref{cycle}), which
can be seen as a vector of the $\mathbb{Q}$-vector space $E=
\mathbb{Q} z_1 \oplus \ldots \oplus \mathbb{Q} z_m$. Note that
given a vector of $E$, we can always get a multiple with integer
coordinates.

As before, let us choose $\tau\in \Gamma_m(f)$ and
choose the numbering of the branches such that $\tau=(1,2,\ldots,m)$.
Then, for every positive integer $l|m$ the candidates to imprimitivity systems of $G_f$ are
\[
\mathcal{B}_l = \{ B_r \}_{r=1, \ldots, d},
\]
where
\[
B_r = \{ r, r+d, \ldots, r + (l-1)d\}
\]
for $d=m/l$. By Proposition~\ref{prop:blocks}, if $\mathcal{B}_l$
is an imprimitivity system, we get a decomposition (up to an
equivalence) $f = f_0 \circ h$. Then, by
Remark~\ref{remark:blocks} the projected cycle $h(C(t))$ is
trivial if and only if $\sum_{j \in B_i} n_j=0$ for
$i=1,\ldots,d$. Similarly, we shall say that $\sum_{j \in B_i}
n_j=0$ for $i=1,\ldots,d$ is the projection of the cycle $C(t)$ by
$\mathcal{B}_l$ even when $\mathcal{B}_l$ is not an imprimitivity
system (note that if $f(z)=z^m$, then $\mathcal{B}_l$ is always an
imprimitivity system). These $d$ equations form a homogeneous
linear system in $m$ unknowns $n_j$. As every column in the
associated matrix has just one $1$ and the rest of the
coefficients zero, while every row has at least one $1$, we can
assure that the rank of the system is $d$. In particular, a
non-trivial (integer) solution of the system is a cycle, since the
sum of all equations is the condition for $C(t)$ being a cycle.

Let us check the dimension of the $\mathbb{Q}$-vector space of the
cycles with trivial projection for every $\mathcal{B}_l$.

\begin{prop}
Under the above conditions, the dimension of the
$\mathbb{Q}$-vector space $T$ of the cycles with trivial
projection for every $\mathcal{B}_l$ is $\phi(m)$, where $\phi(m)$
is Euler's totient function, defined as the number of positive
integers less than $m$ that are relatively prime with $m$.
\end{prop}

\begin{rema}
Let us observe that these cycles with trivial projections are not
necessarily totally unbalanced. To make sure that they are we must
ask for another condition to hold: that they are not balanced.
Then the set of totally unbalanced cycles with trivial projections 
is obtained by substracting some codimension-one 
$\mathbb{Q}$-vector subspaces from 
the $\mathbb{Q}$-vector space of cycles with trivial projections.
\end{rema}

\begin{proof}
First, let us point out that it is enough to check those cycles
with trivial projection for every $\mathcal{B}_p$, where $p$ is a
prime divisor of $m$.

Let $m=p_1^{n_1} \cdot \ldots \cdot p_r^{n_r}$ be the
decomposition of $m$ as product of primes $p_1, \ldots, p_r$. For
every prime $p_i$ we get a homogeneous linear system of
$d_i=m/p_i$ linearly independent equations in $p_i$ unknowns.
Every one of those equations gives a $d_i$-periodic vector $v^k$
of $\mathbb{Q}^m$:
\[
v^k_j = \left\{
\begin{matrix}
1 & \mbox{ if } j \in B_k, \mbox{ i.e., if } j \equiv_{d_i} k
\\
 &
\\
0 & \mbox{ otherwise }
\end{matrix} \right.
\]
So for every prime $p_i |m$ we have $d_i$ linearly independent
vectors of $V_{d_i}$, that is, we get a basis for $V_{d_i}$.

The subspace $T$ of $E$ is the solution of the homogeneous linear
system of $d_1+\ldots+d_r$ equations we get for the primes $p_1,
\ldots, p_r$. So, the codimension of $T$ is the rank of the matrix
associated with this system of equations.

If we consider the equations associated with primes $p_i$ and
$p_j$, we have a set of vectors which consists of a basis for
$V_{d_i}$ and a basis for $V_{d_j}$. Then, we get a system of
generators for $V_{d_i}+V_{d_j}$ and, in consequence, the rank of
the matrix associated with those $d_i+d_j$ equations is the
dimension of $V_{d_i}+V_{d_j}$. That is,
\[
\dim\, (V_{d_i}+V_{d_j}) = \dim V_{d_i} + \dim V_{d_j} - \dim\,
(V_{d_i} \cap V_{d_j}).
\]
Now let us observe that $V_{d_i} \cap V_{d_j} = V_{m/p_ip_j}$,
since $g.c.d.(d_i,d_j)=m/p_ip_j$ because $p_i$ and $p_j$ are
relatively prime. Therefore,
\[
\dim\, (V_{d_i}+V_{d_j}) = d_i + d_j - m/p_ip_j.
\]

In the general case, when we consider all the equations we get:
\begin{equation*}
\begin{split}
{\rm codim}\, T = &  \dim\, (V_{d_1}+ \ldots + V_{d_r}) =
\sum_{i=1}^r \dim V_{d_i} - \sum_{1 \leq i_1 < i_2 \leq r} \dim\,
(V_{d_{i_1}} \cap V_{d_{i_2}}) \\
& + \sum_{1 \leq i_1 < i_2 < i_3
\leq r} \dim\, (V_{d_{i_1}} \cap V_{d_{i_2}} \cap V_{d_{i_3}}) +
\ldots\\
 & + (-1)^{r+1} \dim\, (V_{d_1} \cap \ldots \cap V_{d_r}) \\
= & \sum_{i=1}^r \frac{m}{p_i} - \sum_{1 \leq i_1 < i_2 \leq r}
\frac{m}{p_{i_1}p_{i_2}} + \sum_{1 \leq i_1
< i_2 < i_3 \leq r} \frac{m}{p_{i_1} p_{i_2} p_{i_3}} \\
& + \ldots + (-1)^{r+1} \frac{m}{p_1 \cdot \ldots \cdot p_r} = m
\left( \sum_{i=1}^r \frac{1}{p_i} - \sum_{1 \leq i_1 <
i_2 \leq r} \frac{1}{p_{i_1}p_{i_2}} \right.\\
& \left. + \sum_{1 \leq i_1 < i_2 < i_3 \leq r} \frac{1}{p_{i_1}
p_{i_2} p_{i_3}} + \ldots + (-1)^{r+1}
\frac{1}{p_1 \cdot \ldots \cdot p_r}  \right) \\
= & \, m \left( 1- \prod_{i=1}^r \left( 1 - \frac{1}{p_i} \right)
\right) = m -m \prod_{i=1}^r \left( 1 - \frac{1}{p_i} \right) = m
- \phi(m).
\end{split}
\end{equation*}
Therefore, $T$ is a $\mathbb{Q}$-vector subspace of $E$ of
dimension $\phi(m)$. Consequently, for every $m$ we have $\phi(m)$
(up to a multiple) linearly independent cycles with trivial
projections.
\end{proof}

\begin{exam}\label{cicloTrivializante}
Given $m = 210 = 2 \cdot 3 \cdot 5 \cdot 7$, for the polynomial
$f(z)=z^{210}$ we can find at least $\phi(210)-1=47$ linearly independent
totally unbalanced cycles such that $g(z)= z^2 + z^3 + z^5 + z^7$
is a solution of~(\ref{eq:nullAbel}) but it cannot be reduced to
less than a sum of four reduced solutions. Then, in this case the results
of Pakovich~\cite{PRitt} about the minimum number of reducible
solutions of which the solution of the moment problem is a sum do
not apply.

In general, given $m=p_1 \cdot \ldots \cdot p_s$, for the
polynomial $f(z)= z^m$ we can find at least $\phi(m)-1$ linearly independent
totally unbalanced cycles such that $g(z)= z^{p_1} + z^{p_2} +
\ldots + z^{p_s}$ is a solution of~(\ref{eq:nullAbel}), for $s$
arbitrarily large.
\end{exam}

\section{Balanced Cycles}\label{balanced}

We recall that in order to obtain a complete solution of the tangential zero-dimensional center problem it only remains to solve the problem for balanced cycles. We shall give some examples showing that in this case there exist solutions that are not sums of reduced solutions.

First we recall some definitions and give a characterization of balanced cycles.

Given a field extension $k \to L$ and an algebraic element $\alpha \in L$ over $k$,
the {\em minimal polynomial} of $\alpha$ is the monic polynomial
$p(z) \in k[z]$ of the least degree such that $p(\alpha)=0$. The
minimal polynomial is irreducible over $k$ and any other non-zero
polynomial $q(z) \in k[z]$ with $q(\alpha)=0$ is a multiple of
$p(z)$ in $k[z]$.

The {\em $m$-th roots of unity}, that is, the roots of the polynomial $z^m-1 \in k[z]$ form an abelian cyclic group denoted by $\mu_m$. If the characteristic of $k$ is zero, then $\mu_m$ has order $m$ and
\[
\mu_m = \{ e^{2 \pi i/m} = \epsilon_m, \epsilon_m^2, \ldots, \epsilon_m^m=1 \},
\]
where $\epsilon_m$ is called a {\em primitive} $m$-th root of unity because it generates $\mu_m$. When $k=\mathbb{Q}$, then there are exactly $\phi(m)$ primitive $m$-th roots of unity.

The {\em $m$-th cyclotomic polynomial} is defined to be the minimal polynomial $\Phi_m(z)$ that has the primitive $m$-th roots of unity as simple roots:
\[
\Phi_m(z) =  \underset{\substack{\quad 1\leq k<m,\\gcd(k,m)=1}}{\prod} ( z -
e^{\frac{2 \pi i}{m}k} ).
\]
Obviously, $\phi(m)$ is the degree of $\Phi_m(z)$. Cyclotomic
polynomials are irreducible polynomials with integer coefficients
and $\Phi_m(z)$ is the minimal polynomial of $e^{2 \pi i/m}$ over
$\mathbb{Q}$.

Given a divisor $d$ of $m$, the primitive $d$-th roots of unity
are the elements of $\mu_m$ of order $d$. As the order of any
element of $\mu_m$ is a divisor of $m$, we can conclude that
\[
\mu_m = \underset{d | m}{\coprod}  \{ \text{ primitive $d$-th roots of unity } \}.
\]
Then
\[
z^m - 1 = \underset{\alpha \in \mu_m}{\prod} (z - \alpha) =
\underset{d | m}{\prod} \Phi_d(z).
\]

Given a cycle $C(t)$ by (\ref{cycle}) and a permutation cycle
$\sigma =(p_1,\ldots,p_m) \in \Gamma_m(f)\subset G_f$, we
associate the principal ideal 
\[
I_{C,\sigma}=( [ n_{p_1}z + n_{p_2}z^2 + \ldots + n_{p_m} z^m ] ) \subset
\mathbb{C}[z]/(z^m-1).
\]
The definition above does not depend on the way we represent $\sigma$. Indeed,
if we write $\sigma=(p_2,p_3,\ldots,p_m,p_1)$, then 
\[
\begin{split}
n_{p_1} z + n_{p_2} z^2 + \ldots + n_{p_m} z^m  = & (n_{p_2}z + n_{p_3}z^2+ \ldots +
n_{p_m} z^{m-1} + n_{p_1}z^m)z \\ & - n_{p_1}z(z^m-1),
\end{split}
\]
and consequently
\[
[ n_{p_1} z + n_{p_2} z^2 + \ldots + n_{p_m} z^m  ] =  [ (n_{p_2}z +
n_{p_3}z^2+ \ldots + n_{p_m} z^{m-1} + n_{p_1}z^m) ] [z],
\]
where $[z]$ is invertible in $\C[z]/(z^m-1)$.

We will define $P_{C,\sigma}(z)$ to be any generator of $I_{C,\sigma}$, that is,
\begin{equation}\label{pc}
P_{C,\sigma}(z) = [n_{p_1} + n_{p_2}z + \ldots + n_{p_m} z^{m-1}].
\end{equation}
Note that the cycle $C(t)$ is balanced if and only if $\epsilon_m$ is a root of every polynomial
of $I_{C,\sigma}$ for every $\sigma\in \Gamma_m(f)$, which
is equivalent to $\epsilon_m$ being a root of
$P_{C,\sigma}(z)$ for every $\sigma \in \Gamma_m(f)$.

\begin{prop}\label{prop:balanced}
The cycle $C(t)$ is balanced if and only if
\[
P_{C,\sigma}(z)=\Phi_1(z) \Phi_m(z) d_{\sigma}(z) \quad \forall
\sigma \in \Gamma_m(f),
\]
where $\Phi_m(z)$ is the $m$-th cyclotomic polynomial and
$d_{\sigma}(z) \in \C[z]$ .
\end{prop}

\begin{proof}
Let $C(t)=\sum_{i=1}^m n_i z_i(t)$ be a cycle, $\sigma \in
\Gamma_m(f)$ and $P_{C,\sigma}(z)$ be the polynomial associated
with $C(t)$ and $\sigma$ by~(\ref{pc}). Since $C(t)$ is a cycle,
it follows that $\sum_{i=1}^m n_i=0$, that is, $1$ is a root of
$P_{C,\sigma}(z)$, thus giving that $\Phi_1(z)=z-1$ is a divisor
of $P_{C,\sigma}(z)$.

Now the cycle is balanced if and only if
$P_{C,\sigma}(\epsilon_m)=0$ for every $\sigma \in \Gamma_m(f)$,
that is, if and only if $\epsilon_m$ is a root of the polynomial
$P_{C,\sigma}(z)$. As $\Phi_m(z)$ is the minimal polynomial of
$\epsilon_m$, it is equivalent to $\Phi_m(z)$ dividing
$P_{C,\sigma}(z)$. The claim follows as $\Phi_1(z)$ and
$\Phi_m(z)$ are relatively prime.
\end{proof}

\begin{coro}
If $m$ is prime, then every cycle is totally unbalanced.
\end{coro}

\begin{proof}
First, if $m$ is prime, then $\Phi_m(z)=1+z+\ldots+z^{m-1}$. So
it is not a factor of $P_{C,\sigma}(z)$, and the cycle is
unbalanced. Moreover, since $m$ is prime, then the unique
imprimitivity systems are the trivial ones, and in consequence, $C(t)$ is totally
unbalanced.
\end{proof}

Now we give some examples of solutions when the cycle is balanced.
We will extend them to examples where the cycle is unbalanced but
not totally unbalanced with respect to $f$.

The first example is $f(z)=z^m$, which can be completely solved.
In this case $P_{C,\sigma}(z)$ does not depend on
$\sigma\in\Gamma_m(f)$, so to simplify we write $P_C(z)$.

\begin{prop}\label{prop:xn}
Assume that $f(z)=z^m$ and that $C(t)$ given by \eqref{cycle} is balanced. Then
\[
g(z)=\sum_{j=0}^\ell b_j z^j
\]
is a solution of \eqref{eq:nullAbel} if and only if $b_j=0$
whenever the $\frac{m}{k}$-th cyclotomic polynomial
$\Phi_{m/k}(z)$ does not divide $P_C(z)$, where $k = gcd(m,j)$,
$j>1$.
\end{prop}

\begin{proof}
The $m$-th roots of $f(z)=t$ are $z_i(t) = \epsilon_m^{i-1} s$ for
$i=1,\ldots,m$, where $s$ is an $m$-th root of $t$ and $\epsilon_m$ is a primitive $m$-th root of unity.

Let us compute the value of the expression $\int_{C(t)} g$:
\[
\begin{split}
\int_{C(t)} g &= \sum_{i=1}^m n_i g(z_i) = \sum_{i=1}^m n_i \left( \sum_{j=0}^\ell b_j z_i(t)^j \right) \\
&= \sum_{i=1}^m n_i \left( \sum_{j=0}^\ell b_j ( \epsilon_m^{i-1} s)^j \right) =
\sum_{j=0}^\ell b_j s^j \left( \sum_{i=1}^m n_i \epsilon_m^{(i-1)j} \right)=\\
&=\sum_{j=0}^\ell b_j s^j P_C(\epsilon_m^j).
\end{split}
\]
Hence, $\int_{C(t)}g\equiv 0$ if and only if  $b_j=0$ whenever
$P_C(\epsilon_m^j)\ne0$, $j=0,\ldots,\ell$. For $j=0$ and $j=1$,
$P_C(\epsilon_m^j)=0$ because $C(t)$ is a cycle and because it is balanced respectively.
Then, we can suppose $j>1$. For $j|m$, $\epsilon_m^j$ is a primitive
$\frac{m}{j}$-th root of unity. As we saw in Proposition~\ref{prop:balanced},
$\epsilon_n$ is a root of $P_C(z)$ if and only if $\Phi_n(z)$ divides $P_C(z)$.
Now $P_C(\epsilon_m^j) = P_C(\epsilon_{m/k}^{j/k}) \ne0$ is
equivalent to the condition that the $\frac{m}{k}$-th cyclotomic
polynomial $\Phi_{m/k}(z)$ does not divide $P_C(z)$, where
$k=gcd(m,j)$.
\end{proof}

\begin{exam}\label{deg4}
Let us study some degree four examples. First note that since
$\Phi_1(z)\Phi_4(z)$ $= (z-1)(z^2+1)=z^3-z^2+z-1$ has degree three, by
Proposition \ref{prop:balanced}, $P_C(z)=a(z^3-z^2+z-1)$, $a\ne0$.
Then, the only candidates for balanced cycles are
\[
C(t)=a\left( z_1(t)- z_2(t)+z_3(t)-z_4(t)\right),\quad  a\neq 0.
\]

Take $f(z)=z^4$. Then, by Proposition~\ref{prop:xn}, $g$ is a
solution if and only if it is of the form
\[
g(z)=g_0(z^4)+z g_1(z^4)+z^3 g_3(z^4),
\]
with $g_0,g_1,g_3$ polynomials. Indeed, the only condition
corresponds to coefficients $b_j$, with $j \equiv_4 2$. Then
$k=gcd(4,j)=2$, $\Phi_2(z)=z+1$ does not divide $P_C(z)$ and hence
there are no terms of powers $\equiv_4 2$ in $g$.

This gives examples of solutions $g$ of \eqref{eq:nullAbel} which
can not be written as a sum of $g_k$ such that $f=f_0\circ h$
and $g_k=g_{k,0}\circ h$. For instance, $g(z)=z^3$ is a
solution of \eqref{eq:nullAbel}, while the unique non-trivial
decomposition  of $z^4$ (up to Moebius transformations) is
$z^4=\left(z^2\right)^2$.
\end{exam}

\begin{exam}\label{ex:x4}
Now let $f(z)=T_4(z)=8z^4-8z^2+1$ be the fourth Chebyshev
polynomial and again assume that $C(t)$ is a balanced
cycle, thus,
\[
C(t)=a\left( z_1(t)- z_2(t)+z_3(t)-z_4(t)\right),\quad  a\neq 0.
\]
By direct computation,
\[
\begin{split}
z_1(t)=\sqrt{ 1/2 + \sqrt{1 + t}/(2 \sqrt{2})},\quad 
&z_2(t)=\sqrt{ 1/2 - \sqrt{1 + t}/(2 \sqrt{2})},\\
z_3(t)=-\sqrt{1/2 + \sqrt{1 + t}/(2 \sqrt{2})},\quad 
&z_4(t)=-\sqrt{ 1/2 - \sqrt{1 + t}/(2 \sqrt{2})}.
\end{split}
\]
Rewrite $g=\sum_{i=0}^3 a_i z^{i}g_i(T_4(z))$, for certain
$g_0,g_1,g_2,g_3\in \mathbb{C}[z]$.
Then the Abelian integral along the cycle $C(t)$ reduces to
\[
\sum_{k=1}^4 n_k g(z_k(t))=a_2 a \sqrt{2 + 2 t} g_2(t).
\]
Therefore $g$ is a solution if and only if there exist
polynomials $g_0,g_1,g_3$ such that
\[
g(z)=g_0(T_4(z)) +z g_1(T_4(z))+z^3 g_3(T_4(z)).
\]
Observe that the weak composition conjecture does not hold in this case: take, for instance, the solution $g(z)=z^3$, which does not share any non-trivial compositional factor with $T_4(z) = T_2(T_2(z))$.
\end{exam}

\begin{exam}\label{ex:x6}
Consider now $f(z)=z^6$ and a balanced cycle $C(t)$. Note that
this implies that the associated polynomial $P_C(z)$ must have the
factors $\Phi_1(z)$ and $\Phi_6(z)$. By Proposition~\ref{prop:xn},
the solutions depend on which cyclotomic polynomials divide
$P_C(z)$. The unique candidates are $\Phi_2(z)$ and $\Phi_3(z)$.
Moreover, it is not possible that both $\Phi_2(z)$, $\Phi_3(z)$
divide $P_C(z)$, since $z^6-1 = \Phi_1(z) \Phi_2(z) \Phi_3(z)
\Phi_6(z)$. Therefore there are three possibilities: neither
$\Phi_2(z)$ nor $\Phi_3(z)$ divide $P_C(z)$ (e.g., $C=C_1$),
$\Phi_2(z)$ divide $P_C(z)$ (e.g., $C=C_2$), or $\Phi_3(z)$ divide
$P_C(z)$  (e.g., $C=C_3$).

\small{
\begin{figure}[h]\label{gr:1}
\subfloat[Cycle $C_1(t)$]{
\psfrag{n1}{$-1$}
\psfrag{b1}{$\hspace{-0.2cm}z_1(t)$}
\psfrag{n2}{$2$}
\psfrag{b2}{$\hspace{-0.1cm}z_2(t)$}
\psfrag{n3}{$-2$}
\psfrag{b3}{$z_3(t)$}
\psfrag{n4}{$\hspace{-0.2cm} \phantom{-}1$}
\psfrag{b4}{$z_4(t)$}
\psfrag{n5}{$0$}
\psfrag{b5}{$z_5(t)$}
\psfrag{n6}{$0$}
\psfrag{b6}{$\hspace{-0.1cm}z_6(t)$}
\includegraphics[height=3.0cm]{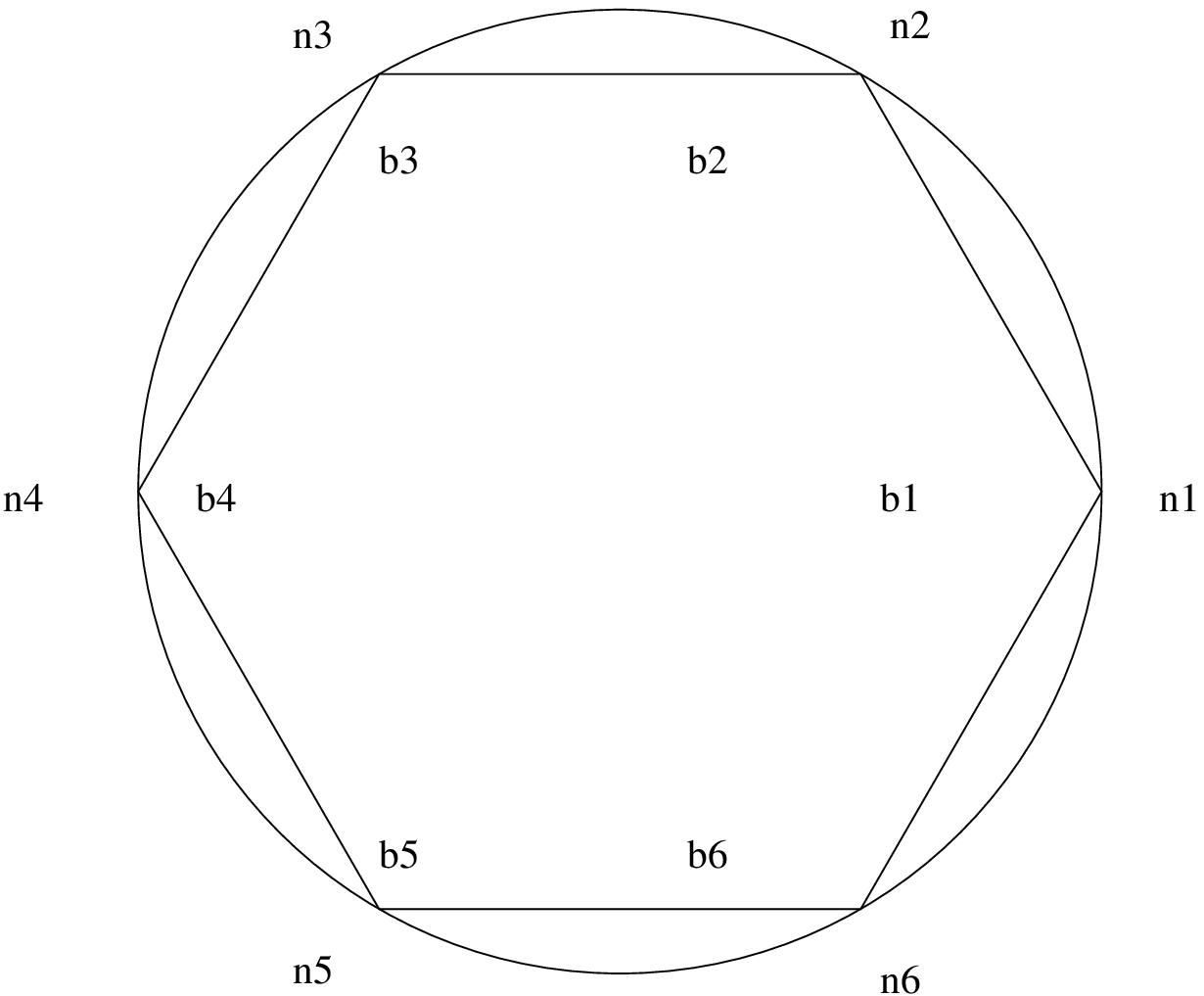}}
$\quad$
\subfloat[Cycle $C_2(t)$]{
\psfrag{n1}{$2$}
\psfrag{b1}{$\hspace{-0.2cm}z_1(t)$}
\psfrag{n2}{$-1$}
\psfrag{b2}{$\hspace{-0.1cm}z_2(t)$}
\psfrag{n3}{$\hspace{-0.1cm}-1$}
\psfrag{b3}{$z_3(t)$}
\psfrag{n4}{$\hspace{0.2cm}2$}
\psfrag{b4}{$z_4(t)$}
\psfrag{n5}{$\hspace{-0.1cm}-1$}
\psfrag{b5}{$z_5(t)$}
\psfrag{n6}{$-1$}
\psfrag{b6}{$\hspace{-0.1cm}z_6(t)$}
\includegraphics[height=3.0cm]{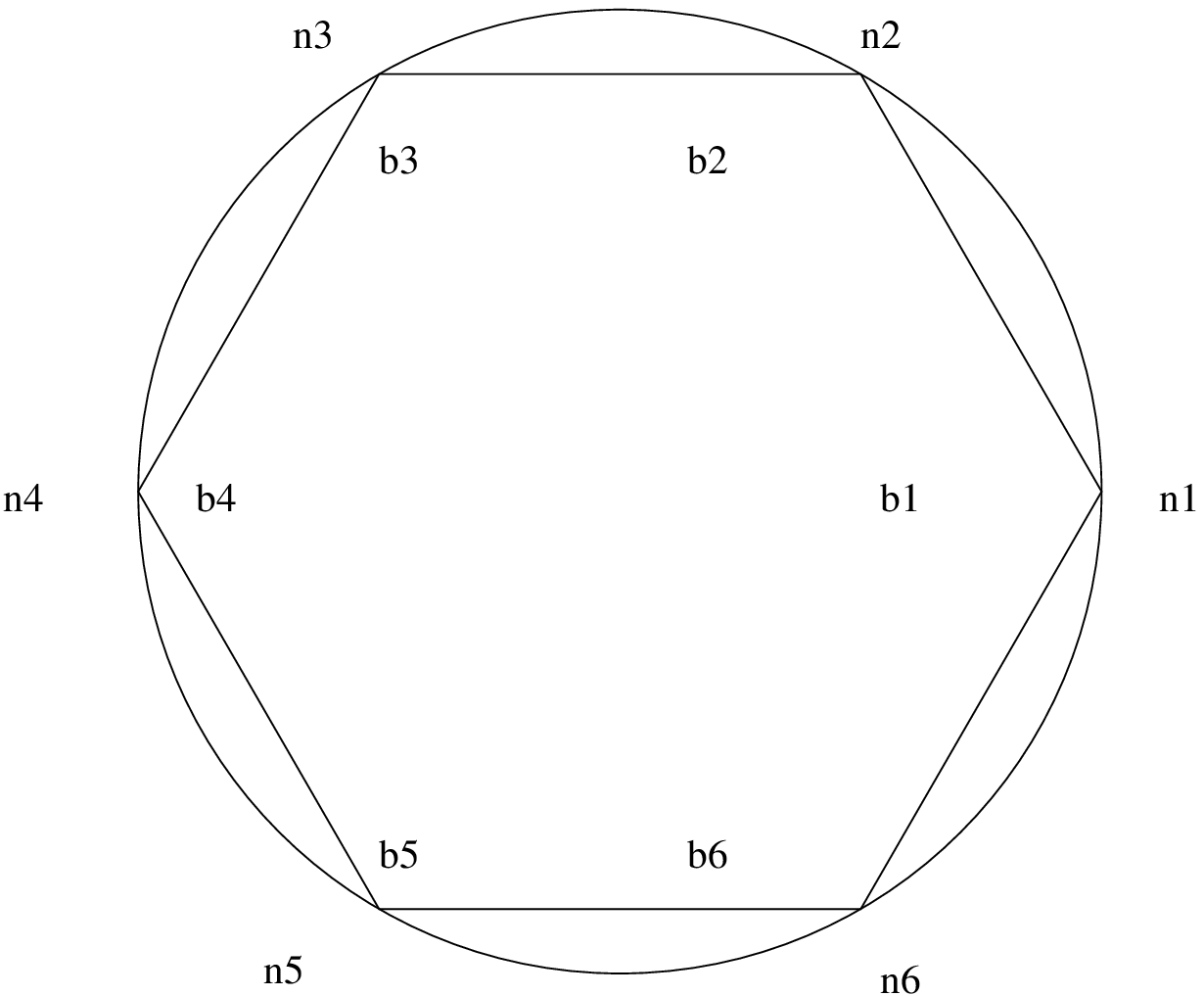}}
$\quad$
\subfloat[Cycle $C_3(t)$]{
\psfrag{n1}{$1$}
\psfrag{b1}{$\hspace{-0.2cm}z_1(t)$}
\psfrag{n2}{$-1$}
\psfrag{b2}{$\hspace{-0.1cm}z_2(t)$}
\psfrag{n3}{$1$}
\psfrag{b3}{$z_3(t)$}
\psfrag{n4}{$\hspace{-0.2cm}-1$}
\psfrag{b4}{$z_4(t)$}
\psfrag{n5}{$1$}
\psfrag{b5}{$z_5(t)$}
\psfrag{n6}{$-1$}
\psfrag{b6}{$\hspace{-0.1cm}z_6(t)$}
\includegraphics[height=3.0cm]{base.eps}}

\caption{Some balanced cycles with six points}
\end{figure}
}

Now, Proposition~\ref{prop:xn} gives us all the solutions:
\begin{enumerate}
\item For the cycle
\[C_1(t)=-z_1(t)+2z_2(t)-2z_3(t)+ z_4(t),\]
$g(z)=\sum_{k\geq0} b_k z^k$ satisfies $\int_{C_1(t)} g\equiv 0$
if and only if $b_k=0$ for every $k\equiv_6 2,3,4$.

\vspace{0.5cm}
\item For the cycle
\[C_2(t)=2 z_1(t)-z_2(t)-z_3(t)+2
z_4(t)-z_5(t)-z_6(t),\] $g(z)=\sum_{k\geq0} b_k z^k$ satisfies
$\int_{C_1(t)} g\equiv 0$ if and only if $b_k=0$ for every
$k\equiv_6 2,4$.
\vspace{0.5cm}
\item For the cycle
\[C_3(t)=z_1(t)-z_2(t)+z_3(t)-z_4(t)+z_5(t)-z_6(t),\]
$g(z)=\sum_{k\geq0} b_k z^k$ satisfies $\int_{C_1(t)} g\equiv 0$
if and only if $b_k=0$ for every $k\equiv_6 3$.
\end{enumerate}
\end{exam}

\begin{exam}
To end with this section, let us show an example of an unbalanced cycle which is
not totally unbalanced.

Take $f(z)=(z^2+az+b)^6$, with $a,b\in\mathbb{C}$.
Let $\tau$ be the element of the monodromy group corresponding
to a simple clockwise loop around all critical values.
We shall number the branches in order to have $\tau=(1,2,\ldots,12)$.

Consider the cycle $C(t)=2 z_1(t)-z_2(t)-z_3(t)+2
z_4(t)-z_5(t)-z_6(t)$. It is easy to check that this cycle is
unbalanced, and that the non-trivial decompositions
of $f$ are $f(z)=z^2\circ (z^2+az+b)^3$,
$f(z)=z^3\circ (z^2+az+b)^2$,
and $f(z)=z^6\circ (z^2+az+b)$.

Let us denote
\[h_1=(z^2+az+b)^3,\quad
h_2=(z^2+az+b)^2,\quad h_3=z^2+az+b.
\]
Then $h_1(C(t))$ is trivial, $h_2(C(t))=4 \bar w_1(t)-2 \bar w_2(t)-2
\bar w_3(t)$ is unbalanced, and $h_3(C(t))=2 w_1(t)-w_2(t)-w_3(t)+2
w_4(t)-w_5(t)-w_6(t)$ is balanced.

By Theorem~\ref{theo:2}, for any $g\in\mathbb{C}[z]$, $\int_C g\equiv 0$ iff there exist
$g_1,g_3\in\mathbb{C}[z]$ such that $g=g_1\circ h_1+g_3\circ h_3$,
and
\[
\int_{h_3(C(t))}g_3\equiv 0.
\]
By Example~\ref{ex:x6}, this is equivalent to
$g_3(z)=\sum_{k\geq0} b_k z^k$ and $b_k=0$ for every $k\equiv_6 2,4$.
\end{exam}

\section{Moment Problem}\label{sec:moment}

Let us consider the following Abel equation
\begin{equation}\label{eq:Abel}
z'=p(w)z^2+\epsilon q(w)z^3,
\end{equation}
where $p$, $q$ are analytic functions. We say that a solution $z(w)$ of \eqref{eq:Abel} is closed if $z(0)=z(1)$. Let us denote
\[
f(w)=\int_0^w p(t)\,dt,\quad g(w)=\int_0^w q(t)\, dt.
\]

Assume that we have a center at the origin for $\epsilon=0$ (thus, every bounded solution of the system is closed), which happens iff $f(0)=f(1)$, and we want to study when the center persists after perturbation (i.e., for every $\epsilon$ we have a center). This problem comes from the center-focus problem for certain planar systems, which after a change to polar coordinates become of the form of \eqref{eq:Abel} with $p,q$ trigonometric polynomials. The case when $p,q\in\mathbb{C}[w]$  has been extensively studied and
solved. We show in this section how both cases can be seen as zero-dimensional tangential center problems.

To first order in $\epsilon$, the center is persistent if and only if $(\partial z/\partial \epsilon)(1,z_0,0)=0$ for every $z_0$ close to zero, where $z(t,z_0,\epsilon)$ is the solution of \eqref{eq:Abel} determined by $z(0,z_0,\epsilon)=z_0$.

Differentiating $z^{-1}(t,z_0,\epsilon)$ with respect to $\epsilon$ and evaluating at $t=1$, $\epsilon=0$, one has that $(\partial z/\partial \epsilon)(1,z_0,0)=0$ for every $z_0$ close to zero if and only if
\begin{equation}\label{eq:momentproblem}
0=\int_0^{1} \frac{q(w)dw}{1-z_0 f(w)}=\sum_{k=0}^\infty\left(\int_0^{1} f^k(w)q(w)dw\right)z_0^k,
\end{equation}
for every $z_0$ close to zero (see \cite{BFY}). 

The moment problem~\cite{BFY1} asks for a given $p$, to find all $q$ such that \eqref{eq:momentproblem} is satisfied. It has been extensively studied in recent years (see, e.g., \cite{BRY}, \cite{FPYZ} and references therein).

The first case to be studied is the case when $p$ and $q$ are polynomials, recently solved by Pakovich and Muzychuk~\cite{PM}. More precisely,  they prove:

\begin{theo}[\cite{PM}]\label{PM}
Take a polynomial $p$ such that $f(0)=f(1)$. A polynomial $q$ is a solution of the polynomial moment problem if and only if there exist polynomials $f_1,\ldots,f_r,g_1,\ldots,g_r,h_1,\ldots h_r$ such that $h_k(0)=h_k(1)$, $f=f_k\circ h_k$ for every $k$ and
\[
g=\sum_{k=1}^r g_k\circ h_k.
\]
\end{theo}

We shall show that the polynomial moment problem is equivalent to a zero-dimensional tangential center problem for $f$ and $g$ defined as in Theorem~\ref{PM}, and for the cycle
\[
C_f(t)=\sum_{i=1}^{n_0} n_1 z_{a_i}(t)-\sum_{i=1}^{n_1} n_0 z_{b_i}(t),
\]
where $z_{a_i}(t)$ are all the solutions of $f(z_{a_i}(t))=t$  close to $0$ for $t$ close to $f(0)$, and analogously for $z_{b_i}(t)$ and $1$. We shall prove that $C_f(t)$ is totally unbalanced, so Theorem~\ref{PM} follows from Theorem~\ref{theo:2}. Note however that in the proof of Theorem~\ref{theo:2} we use methods developed in Theorem~\ref{PM}.

\begin{prop}\label{prop:trigmoment}
Given polynomials $p,q$ such that $f(0)=f(1)$, the polynomial $q$ is a solution of \eqref{eq:momentproblem} if
and only if
\[
\int_{C_f(t)} g \equiv 0.
\]
Moreover, $C_f(t)$ is totally unbalanced.
\end{prop}

\begin{rema}
As a consequence of Proposition~\ref{prop:trigmoment}, if $\int_{C_f(t)} g \equiv 0$, then $g$ is totally
determined by Theorem~\ref{theo:2}.
\end{rema}

\begin{proof}
By Theorem~1.1 of \cite{PM}, $q$ is a solution of the polynomial moment problem if and only if there exists polynomials $f_1,\ldots,f_r,g_1,\ldots,g_r,h_1,\ldots h_r$ such that $h_k(0)=h_k(1)$, $f=f_k\circ h_k$ for every $k$ and
\[
g=\sum_{k=1}^r g_k\circ h_k.
\]

On the other hand, if $C_f(t)$ is totally unbalanced, then Theorem~\ref{theo:2}
implies that $g$ is a solution of \eqref{eq:momentproblem}
for the cycle $C_f(t)$ if and only if there exist polynomials
$f_1,\ldots,f_s,g_1,\ldots,g_s,h_1,\ldots h_s$
such that $f=f_k\circ h_k$ for every $k$,
\[
g=\sum_{k=1}^r g_k\circ h_k,
\]
and $h_k(C_f(t))$ is trivial.

First, let us prove that $h_k(0)=h_k(1)$ is equivalent to $h_k(C_f(t))$ being trivial. To this end,
let $\mathcal{B}$ be the imprimitivity system corresponding to some
$h_k$ such that $f=f_k\circ h_k$.

If $h_k(0)=h_k(1)$, then the blocks of $\mathcal{B}$ contain
elements of both $\{a_i\}$ and $\{b_i\}$: take $t$ a non-critical
value of $h_k$ and joint it to $h_k(0)=h_k(1)$ by an arc. Then
$h_k^{-1}$ maps this arc into some arcs close to $0$ and others
close to $1$. Let $\sigma\in G_f$ denote an element of the
monodromy group corresponding to a cycle around
$h_k(0)=h_k(1)$. Reorder $a_i,b_i$ so that $\sigma$ permutes
cyclically $a_1,\ldots,a_{n_0}$ and $b_1,\ldots,b_{n_1}$ and
that $a_1,b_1$ belongs to the same block of $\mathcal{B}$. Now,
the blocks reduced to $a_1,\ldots,a_{n_0}$, $b_1,\ldots,b_{n_1}$
consist of the congruence classes modulo $d$. Therefore, for every
block $B\in \mathcal{B}$,
\[
\sum_{k\in B}n_k=\sum_{i=1}^{n_0/d}n_1-\sum_{i=1}^{n_1/d}n_0=0.
\]

Conversely, if $h_k(0)\neq h_k(1)$, since the blocks contain
branches where the value of $h_k$ is the same, then every block
$B$ of the imprimitivity system $\mathcal{B}$ corresponding to
$h_k$ either contains elements of $\{a_i\}$ or elements of
$\{b_i\}$, and then $\sum_{k\in B}n_k\neq 0$, and $h_k(C_f(t))$ is
not trivial.

Therefore, to conclude it only remains to prove that $C_f(t)$ is totally unbalanced.
Indeed, let us choose $\tau$ a permutation corresponding to a simple loop around infinity and 
reorder the numbering of the branches such that $\tau=(1,2,\ldots,m)$. Define
\[
V(0)=\{\epsilon_m^{a_1},\ldots,\epsilon_m^{a_{n_0}}\},\quad
V(1)=\{\epsilon_m^{b_1},\ldots,\epsilon_m^{b_{n_1}}\}.
\]
The Monodromy Lemma in~\cite{P2} states that the convex hulls of the sets $V(0)$ and
$V(1)$ are disjointed when $f(0)=f(1)$, thus, $\mathbb{C}$ can be
divided into two half-planes $P^1,P^2$ such that $V(0)\in P^1$ and
$V(1)\in P^2$.

Now, the center of mass of $C_f(t)$ is the origin (i.e.,
$P_{C_f,\tau}(\epsilon_m)=0$) if and only if the center of mass of
$V(0)$ is equal to the center of mass of $V(1)$. But the center
of mass of $V(0)$ belongs to $P^1$ and the center of mass of
$V(1)$ belongs to $P^2$. Therefore, $C_f(t)$ is unbalanced.

If $f=f_k\circ h_k$, then there are two possibilities for the
projected cycle $h_k(C_f(t))$: If $h_k(0)=h_k(1)$, then
$h_k(C_f(t))$ is trivial. If $h_k(0)\neq h_k(1)$, then we shall
prove that $h_k(C_f(t))$ is unbalanced. Indeed, the permutation
induced on the branches $z_{a_1},z_{a_2},\ldots,$ $z_{a_{n_0}}$ by
the preimage of a loop around $f(0)$ is a permutation cycle.
Therefore, the permutation induced on the branches
$h_k(z_{a_1}),h_k(z_{a_2}),\ldots h_k(z_{a_{n_0}})$ must also be a
permutation cycle. Let us denote by $d_a$ the order of that
permutation cycle and by $d_b$ the order of the permutation cycle
induced by a loop around $f(1)$ on
$h_k(z_{b_1}),h_k(z_{b_2}),\ldots h_k(z_{b_{n_1}})$. Then
\[
h_k(C_f(t))=\sum_{j=1}^{n_0/d_a} d_a n_1 h_k(z_{a_j}(t))-\sum_{j=1}^{n_1/d_b} d_b n_0  h_k(z_{b_j}(t)).
\]
Applying the same arguments as for $C_f(t)$, one obtains that $h_k(C_f(t))$ is unbalanced.
As a consequence, $C_f(t)$ is totally
unbalanced with respect to $f$.
\end{proof}

Now, consider the trigonometric version of the moment problem. Given a
trigonometric polynomial $f$, find all trigonometric polynomials $g$ such that
\begin{equation}\label{eq:trigmomentproblem}
\int_0^{2\pi} f^k(w)g'(w)dw=0,\quad k\geq 0,
\end{equation}

By the change of variable $z\to e^{iz}$, the trigonometric moment
problem is equivalent to the following problem (see \cite{ABC}):
Given a Laurent polynomial $f$ ($f\in \mathbb{C}[z,z^{-1}]$),
find all $g\in \mathbb{C}[z,z^{-1}]$ such that
\begin{equation}\label{eq:laurentproblem}
\int_{|z|=1} f^k(w)g'(w)dw=0,\quad k\geq 0.
\end{equation}
Moreover, if $\int_{|z|=1} g'(w)\,dw=0$,
then by integration by parts,
\[
\int_{|z|=1} f^k(w)g'(w)dw=-k\int_{|z|=1} f^{k-1}(w)f'(w)g(w)dw,\quad k\geq 1.
\] 
Since $\int_{|z|=1} f^k(w)f'(w)g(w)dw$, $k \geq 0$,
are the coefficients of the power series at $t=\infty$ of
\[
I(t)=\int_{|z|=1} \frac{g(w)f'(w)\,dw}{t-f(w)}.
\]
then for any fixed Laurent polynomial $f$, a Laurent polynomial $g$
such that  $\int_{|z|=1} g'(w)\,dw=0$
 is a solution of \eqref{eq:laurentproblem}
if and only if $I(t)\equiv 0$. 

\begin{prop}
Let
\[
f(z)=\sum_{k=-n}^m c_k z^k
\]
be a proper Laurent polynomial, that is, $n,m\geq 1$. A Laurent polynomial $g$ 
such that $\int_{|z|=1} g'(w)\,dw=0$
is a solution of \eqref{eq:laurentproblem} if and only if
$g$ is a solution of \eqref{eq:nullAbel} for the cycle
\[
C_f(t)=\sum_{k=1}^n m z_{a_k}(t)-\sum_{k=n+1}^{n+m} n z_{b_k}(t),
\]
where $z_{a_i}(t)$ (resp. $z_{b_i}(t)$)  are the branches close to $0$ (resp. $\infty$) 
of a $t$ close to infinity.
\end{prop}

\begin{proof}
Take $t$ close to infinity. By the Residue Theorem,
\[
I(t)=2\pi i\left(-\sum_{k=1}^n g(z_{a_k}(t))+Res(gf'/(t-f),0))\right).
\]
Assume that $I(t)\equiv 0$ and take $\sigma\in G_f$. By analytic continuation,
\[
\sum_{k=1}^n g(z_{a_k}(t))-\sum_{k=1}^n g(z_{\sigma(a_k)}(t))=0
\]
Since $G_f$ is transitive, if we sum the previous formula over $G_f$, we obtain
\[
\sum_{k=1}^n |G_f| g(z_{a_k}(t))-\sum_{k=1}^{n+m}\frac{|G_f|n}{n+m} g(z_{k}(t))=0.
\]
Dividing by $|G_f|$ and multiplying by $n+m$ gives
\[
\int_{C_f(t)} g=\sum_{k=1}^n m g(z_{a_k}(t))-\sum_{k=n+11}^{n+m} n g(z_{b_k}(t))=0.
\]

Conversely, if $g$ is a solution of \eqref{eq:nullAbel}, then
using the Residues Theorem in the two regions of the complementary of the unit
circle, 
\[
I(t)=\frac{2\pi i}{n+m}(m Res(gf'/(t-f),0))-n Res(gf'/(t-f),\infty))).
\]
Now, take $\alpha=1-(t z^n -f(z)z^n)/c_n$. Note that
$\alpha\to 0$ when $z\to 0$. Then for $z$ close to $0$,
\[
\frac{1}{t-f(z)}=\frac{z^n}{c_n(1-\alpha)}=\frac{z^n}{c_n}\sum_{k=0}^\infty \alpha^k
=\frac{z^n}{c_n}\sum_{k=0}^\infty \left(1-\frac{t z^n -f(z)z^n}{c_n} \right)^k.
\]
Therefore, $Res(gf'/(t-f),0)$ is a polynomial in $t$, and the same holds for $Res(gf'/(t-f),\infty)$. 
To conclude, note that $I(t)\to 0$ as $t\to \infty$,
then $I(t)\equiv 0$.
\end{proof}

In the general rational case, 
Pakovich~\cite{Prat} has proved that if $\gamma$ is a rational curve and $f,g$ are rational functions, 
then
\[
\int_{\gamma} f^k(z)\,dg(z)=0,\quad \text{for every }k\geq 0,
\]
if and only if $\int_{C_{f,k}(t)} g\equiv 0$ for a finite number
of cycles $C_{f,k}(t)$ defined in terms of the the ``dessin d'enfants'' of $f$.

\section{The Tangential Center-Focus Problem in the Hyper-elliptic Case}\label{hyper}

Let us now consider the original center problem in the two-dimensional
space. Let $F(x,y)\in\mathbb{C}[x,y]$ and consider the foliation
\begin{equation}\label{pertF}
dF+\epsilon\omega=0
\end{equation}
deforming the initial foliation $F=t$, where $\omega$ is a one-form. Let  $t_0$ be a regular
value of $F$ and let $\gamma(t_0)\subset F^{-1}(t_0)$ be a closed
path in the leaf $F^{-1}(t_0)$. Take a transversal $T$ to the
leaves of $F$ parameterized by values $t$ of $F$ and consider the
holonomy along $\gamma(t_0)$ and the corresponding displacement
map $\Delta(t,\epsilon)$ (holonomy minus identity). As the path
$\gamma(t_0)$ pushes to nearby closed paths $\gamma(t)\subset
F^{-1}(t)$, we have that $\Delta(t,0)\equiv 0$. We say that the
family $\gamma(t)$ is a center. We have
\begin{equation}\label{Delta1}
\Delta(t,\epsilon)=-\epsilon\int_{\gamma(t)}\omega+o(\epsilon).
\end{equation}
Of course the Abelian integral in \eqref{Delta1} depends only on
the homology class of $\gamma(t)\in H_1(F^{-1}(t))$. We say that the family
$\gamma(t)$ is a persistent center if
$\Delta(t,\epsilon)\equiv0$ and a tangential center if
$\int_{\gamma(t)}\omega\equiv0$.

We consider more precisely the \emph{hyper-elliptic} case
$F(x,y)=y^2+f(x)$, $f(x)\in\mathbb{C}(x),$ $\deg(f)=m$. Recall
that the tangential center problem in the hyper-elliptic case for
vanishing cycles was solved in \cite{CM}. Vanishing cycles are a
particular class of simple cycles. We want to study the general
tangential center problem in the hyper-elliptic case.

It is well-known (see for instance \cite{Z}) that any one-form
$\omega$ is relatively cohomologous to the form
\begin{equation}
\omega=g(x,y)dF+dR(x,y)+\kappa(x)ydx,
\end{equation}
for some $g(x,y), R(x,y)\in\mathbb{C}[x,y]$, $\kappa(x)\in\mathbb{C}[x]$. Hence, the tangential
center problem reduces to the problem of characterizing the vanishing of the integral
\begin{equation}
\int_{\gamma(t)}\kappa(x)ydx
\end{equation}
along some cycle $\gamma(t)\in H_1(F^{-1}(t)).$

We represent $F^{-1}(t_0)$ as a Riemann surface $y=\sqrt{t_0-f(x)}$
and consider its first homology group $H_1(F^{-1}(t_0))$. Let
$x_1(t_0),\ldots,x_m(t_0)$ be all the roots of $f^{-1}(t_0)$. It
is well-known that the homology of the  fiber ${F^{-1}(t_0)}$ is
generated by \emph{simple cycles} $C_{ij}(t_0)$ going from
$x_i(t_0)$ to $x_j(t_0)$ on one leaf of the Riemann surface
followed by the lift of the same path on the other leaf travelled
in the opposite direction. As the value of $y$ on the second half
of the path is opposite to the value on the first part as well as
the sense of travel, the two halves of the integral add up. Set
\begin{equation}\label{g}
g(x,t)=\int \kappa(x)\sqrt{t-f(x)}dx
\end{equation}
the indefinite integral of $\kappa(x)\sqrt{t-f(x)}$. Up to the
problem of multivaluedness of $g$, this gives that the original
tangential center problem in the plane reduces to the
zero-dimensional tangential center problem
\begin{equation}
\int_{C(t)}g(x,t)=\sum_{i=1}^m n_i g(x_i(t),t)\equiv 0,
\end{equation}
where the zero-cycle $C(t)=\sum_{i=1}^m n_i x_i(t)$ is composed
of the ramification points $x_i(t)\in f^{-1}(t)$ in the
presentation of the class of  $\gamma(t)$ in $H_1(f^{-1}(t))$.
Note that $\sum_{i=1}^m n_i=0$, that is, $C(t)$ is indeed a cycle
as the class of $\gamma(t)$ is generated by simple cycles.

The explicit dependence of $g(x(t),t)$ on $t$ is not a problem. We
must however be more careful here as the function $g(x,t)$ is
multivalued in $x$ due to the multivaluedness of the square root.
It has the same (two-sheeted) Riemann surface as the function
$y=\sqrt{t-f(x)}.$ We will now be more precise in the choices we
make in representing the cycle $\gamma(t)$. We take $t_0$ big,
order cyclically the roots $x_i(t_0)$ and suppose for simplicity
that $m$ is even. We present the Riemann surface of $g$ as a two-sheeted 
surface with cuts joining $x_1(t_0)$ to $x_2(t_0)$, next
$x_3(t_0)$ to $x_4(t_0)$ etc. When crossing any cut a path changes
the sheet. Let $\gamma_{ij}(t_0)$ be the cycle going first from
$x_i(t_0)$ to $x_j(t_0)$ in the upper sheet and then along the same
projection from $x_j(t_0)$ to $x_i(t_0)$ on the lower sheet. We
suppose that the cycle avoids cuts from exterior  (that is, it belongs
to the complement of the convex hull  of the cuts except for the
ramification points $x_i(t_0)$).

If $m$ is odd, we have to add a cut from $x_m$ to infinity, but we
will stick rather to the case when $m$ is even for simplicity. The
space $\mathcal{B}=\{\gamma_{i,i+1}(t_0): i=1,\ldots,m-1\}$
provides a basis of $H_1(F^{-1}(t_0))$.  Any cycle $\gamma(t_0)$
can hence be written in a unique way as
\begin{equation}\label{gamma}
\gamma(t_0)=\sum_{i=1}^{m-1} n_{i}\gamma_{i,i+1}(t_0).
\end{equation}
We associate the zero-cycle
\begin{equation}\label{C}
\phi(\gamma(t_0))=C(t_0)=\sum_{i=1}^{m-1} n_{i} (x_{i+1}(t_0) - x_i(t_0))
\end{equation}
to the one-cycle $\gamma(t_0)$. 
In fact, we have constructed an isomorphism
 \begin{equation}\label{phi}
\phi:H_1(F^{-1}(t_0))\to \tilde H_0(f^{-1}(t_0)).
\end{equation}

Then
\begin{equation}\label{int}
\int_{\gamma(t_0)}\kappa(x)ydx=2\int_{C(t_0) }g.
\end{equation}

We say that a cycle $\gamma(t)$ in $F^{-1}(t)$ is balanced
(unbalanced, totally unbalanced) if the corresponding zero-cycle $\phi(\gamma(t))$ in
$f^{-1}(t)$ is balanced (unbalanced, totally unbalanced).

Let $G_F$ denote the usual action of the monodromy group of $F$ on $H_1(F^{-1}(t_0))$
and let $G_F^0$ denote the conjugated group: \[G^0_F=\{\phi\circ\sigma\circ\phi^{-1}:\sigma\in G_F\}.\] 
Then  $G^0_F$ acts  on $\tilde H_0(f^{-1}(t_0))$.

The action of the group $G^0_F$ is very much related to the action
of the monodromy group $G_f$ of $f$, but in general is more
complicated due to cuts in the Riemann surface of $g$. 

We believe that the following conjecture, analogous to Theorem
\ref{theo:2}, holds:
\begin{conj}\label{conj:1}
Let $F(x,y)=y^2+f(x)\in\mathbb{C}[x,y]$, $\kappa\in\mathbb{C}[x]$ and let $g(x)$ be given by \eqref{g}.

\begin{enumerate}
\item If $\gamma(t)$ is a totally unbalanced cycle of $F$, then
\begin{equation}\label{eq:nullAbel1}
\int_{\gamma(t)} \kappa(x)ydx\equiv 0
\end{equation}
if and only if there exist
$f_1,\ldots,f_s,h_1,\ldots,h_s\in\mathbb{C}[x]$ and analytic
functions in a neighborhood of infinity $g_1,\ldots,g_s$ such that $f=f_k\circ h_k$, $g=g_1\circ
h_1+\ldots+g_s\circ h_s$ and  the cycles $H_k(\gamma)$ are trivial
for every $k=1,\ldots,s$, where $H_k(x,y)=(h_k(x),y).$

\item If $\gamma(t)$ is unbalanced, then
\eqref{eq:nullAbel1}  holds if and only if there exist
$f_1,\ldots,f_s,h_1,\ldots,$ $h_s\in \mathbb{C}[x]$, $g_1,\ldots,g_s$
analytic in a neighborhood of infinity 
such that $f=f_k\circ h_k$, $g=g_1\circ
h_1+\ldots+g_s\circ h_s$, and for every $k$, either the projected
cycle ${H_k}({\gamma})$ is trivial or it is balanced and
\[
\int_{{h_k}({\phi(\gamma(t))})} g_k \equiv 0.
\]
\end{enumerate}
\end{conj}

The proof should go similarly as the proof of Theorem \ref{theo:2}.
The delicate point is the use of Lemma \ref{lem:irreducible}. If in 
Lemma \ref{lem:irreducible} we could replace $G_f$-irreducible and 
$G_f$-invariant by $G_F^0$-irreducible and $G_F^0$-invariant, then we 
could prove the conjecture by using the isomorphism $\phi$ and then following the
proof of Theorem \ref{theo:2}. Here it is essential that the
function $f$ is polynomial. The inductive argument is on the
divisors of the degree of $f$. The fact that the function $g$ is
only analytic is not a problem. The difficulty is that $G_F^0$ is
not a permutation of the roots as $G_f$. In fact it acts on cycles
and not on individual roots.

\begin{exam} 
Let $F(x,y)=y^2+x^m$ and assume that $m$
is even. Let $\gamma(t)$ be a balanced cycle of $F$. Let
$P_{\gamma}$ be the polynomial associated to the cycle $\gamma$.
Note that here the group $G_F^0$ coincides with $G_f$ because
there is only one critical value and hence only one generator
(corresponding to winding around zero). Hence Proposition
\ref{prop:xn} applies and $\int_{\gamma(t)}\kappa(x)ydx$ vanishes
if and only if $b_j(t)=0$, whenever $\Phi_{m/k}(x)$ does not
divide $P_{\gamma}(x)$, where $k=gcd(m,j)$ and $g(x,t)=\sum b_j(t)
x^j$ is the development of $g$ given in \eqref{g} in a
neighborhood of infinity. Following Example \ref{deg4}, for
$F(x,y)=y^2+x^4$, and the balanced cycle
$\gamma(t)=\phi^{-1}(C(t))$, with
$C(t)=a(x_1(t)-x_2(t)+x_3(t)-x_4(t))$, we have that
$\int_{\gamma(t)}\kappa(x)ydx\equiv0$ if and only if the function
$g(x,t)$ given by \eqref{g} is of the form
$g(x,t)=g_0(x^4,t)+xg_1(x^4,t)+x^3g_3(x^4,t)$, where the functions
$g_i$, for fixed $t$ are analytic in a neighborhood of infinity in
$x$.

Recall the definition of $g$.
We use the development in a neighborhood of infinity. 
There $\sqrt{t-x^4}$ develops as $x^2$ times a function of $(x^4,t)$. This gives that the powers of $x$ appearing in $\kappa(x)\sqrt{t-x^4}$ are by two higher than the powers in $\kappa(x)$. One extra shift in powers comes from integration. This gives all together that the powers appearing in $\kappa(x)$ are by $3$ modulo $4$ higher than the powers in $g$.  
Finally we get that
$\int_{\gamma(t)}\kappa(x)ydx\equiv0$ if and only if the
polynomial $\kappa(x)$ is of the form
$\kappa(x)=\kappa_1(x^4)+x^2\kappa_2(x^4)+x^3\kappa_3(x^4)$, with
$\kappa_i\in\mathbb{C}[x]$.
\end{exam}

\section{Tangential canard centers for generalized Van der Pol's equations}\label{sec:sf}

Consider the singular perturbation of the generalized Van der Pol's equation
\begin{equation}\label{eq:vanderpol}
\epsilon x''+\alpha(x)x'+\beta(x)=0,
\end{equation}
for certain rational functions $\alpha,\beta$. 

As in \cite{DR0} we blow up the singular point in the family and study the 
blown up system. A natural question
is when after perturbing with
the blow up parameter $u$ we obtain a center. The tangential version of this
problem is to study when the perturbed equation has
a center ``to first order'' in the parameter $u$.

In this section on the example of Van der Pol's singular
perturbation, we show the relevance of the zero-dimensional
tangential center problem for studying persistence of centers in
slow-fast systems.

First, following~\cite{DR0}, \eqref{eq:vanderpol} can be
transformed into the slow-fast system
\[
\begin{split}
x'&=y-f(x),\\
y'&=\epsilon G(x),
\end{split}
\]
where $f,G$ are rational functions and $\epsilon \in(0,+\infty)$.
We shall assume in addition that $f$ is a Morse function
$f(x)=\frac{1}{2} x^2+o(x^2)$. Note that for every $t>0$ close
enough to zero, $f(x)=t$ has exactly two solutions also close to
zero. Let $z_1(t)<0<z_2(t)$ denote those two solutions. Also, assume
that
\[
G(x)=a_0-x(1+O(x)),
\]
with $1+O(x)>0$ for $x$ close to zero.

Now the system rewrites as
\begin{equation}\label{eq:slowfast}
\begin{split}
x'&=y-f(x),\\
y'&=\epsilon(a_0-x\bar G(x)).
\end{split}
\end{equation}

We blow-up the family. Here we consider only one chart:
\begin{equation}\label{blowup}
(x,y,a_0,\epsilon)=(u\bar x,u^2\bar y,uA_0,u^2).
\end{equation}
Dividing by $u$, the blown up vector field, has a center at $(\bar x,\bar y)=(0,0)$ on the divisor $u=0$, for $A_0=0$. We are interested in the persistence (or rather first order persistence of this center).
With this notation, the second part of Theorem~1 of \cite{DR}
can be formulated:

\begin{prop}[\cite{DR}]\label{prop:DR}
Under above assumptions, for $(A_0,u)$ close to $(0,0)$
$u>0$ there is a $\mathcal{C}^\infty$-function 
$\Delta(\bar y,A_0,u)$ defined on a transversal $\bar x=0$ such 
that the closed orbits are given by $\Delta(\bar y,A_0,u)=0$.

Moreover,
\[
\frac{\partial \Delta}{\partial \bar y}=-I(\bar y)+O(u),
\]
and
\begin{equation}\label{eq:integralslowfast}
I(t)=\int_{z_1(t)}^{z_2(t)} \bar xw(\bar x) F^2(\bar x) \,d\bar x,
\end{equation}
where
\[
w(\bar x)=-\dfrac{\bar x}{\bar G(\bar x)},\quad F(\bar x)=\frac{1}{\bar x}\frac{\partial f}{\partial \bar x}(\bar x).
\]
\end{prop}

If $I(t)\equiv 0$ for some $w(x)$ and all $t>0$, then
we have that \eqref{eq:slowfast} has a center to first order in $u$.
Rewriting \eqref{eq:integralslowfast} as
\[
I(t)=g(z_1(t))-g(z_2(t)),
\]
where $g$
is a primitive of $\bar xw(\bar x)F^2(\bar x)$,
the problem of characterizing $I(t)\equiv 0$ is equivalent to
the zero-dimensional tangential center problem for the
cycle $C(t)=z_1(t)-z_2(t)$.

Proposition~\ref{prop:simplecycle}
states that $g$ is a solution iff there exist $f_0,g_0,h$ with $h(z_1(t))=h(z_2(t))$
such that $f=f_0\circ h$, and $g=g_0\circ h$.
Therefore,
\[
g_0'(h(\bar x))h'(\bar x)=\bar x w(\bar x) F^2(\bar x).
\]
Now, replacing we get
\[
\bar G(\bar x)=-\frac{ \left(\frac{\partial f}{\partial \bar x}(\bar x)\right)^2}{g_0'(h(\bar x))h'(\bar x)}
\]
for any $g_0'$ and any $h$ such that $h(z_1(t,c))=h(z_2(t,c))$ and $f=f_0\circ h$.
Thus, given $f$, it is possible to obtain all $\bar G\in\mathbb{C}(z)$ 
such that the perturbed system has a center at first order in $u$.

\section{Concluding remarks and open problems}

The zero-dimensional tangential center problem (Problem~\ref{problem2}) is a kind of toy
example for studying the tangential center problem for the two-dimensional 
phase space (vanishing of Abelian integrals on one-dimensional cycles). In Section \ref{hyper} we see how the initial problem in the hyper-elliptic case reduces to a problem in
the zero-dimensional case. A first problem is to generalize our
results and prove an analogue of Theorem \ref{theo:2} for general
cycles on hyper-elliptic systems.

To completely solve the zero-dimensional tangential center problem
it remains to solve the problem for balanced cycles. We believe
that the existence of powers $f_i(x)=x^i$ or Chebyshev polynomials
$T_i(x)=\cos(i\arccos(x))$ as composition factors of $f$ plays a center role in balanced tangential centers. The study of balanced centers of $f_i(x)=x^i$ was easy. The Chebyshev
case should be only slightly more complicated. Example~\ref{deg4}
shows that even the weak composition conjecture is not valid for
the zero-dimensional tangential center problem. It would be very
interesting to formulate a new plausible conjecture for necessary
and sufficient condition of a tangential center.

Recall that Abelian integrals represent the linear term of the
displacement function. If they vanish, one searches for the
first non-zero term in the development of the displacement
function. In the one-dimensional cycle case it is known that this
first term is given by iterated integrals (see \cite{Ch},
\cite{G}, \cite{JMP1} and \cite{JMP2}). It would be interesting to
define and study iterated integrals in the zero-dimensional case,
too. It seems that they cannot exist for unbalanced cycles, but
probably for balanced cycles it  makes sense. 

When studying the center problem one often develops the
displacement function at the center. The center condition is the
condition of vanishing of all coefficients in this development.
This gives a growing sequence of ideals which stabilize.
The order at which they stabilize is called the Bautin index. Bautin index was calculated by Bautin for quadratic vector fields \cite{B}. Its calculation for general vector fields is an important unsolved problem. Zero-dimensional problem being more accessible could help to develop techniques for studying this
problem.

Here we solved  the polynomial tangential center problem for
totally unbalanced cycles. It would be interesting to solve it for
rational or Laurent systems. The Laurent systems case corresponds
to trigonometric polynomials.

\end{document}